\documentclass{article}
\usepackage[spanish]{babel}
\usepackage{amsmath,amssymb,graphicx,textcomp}

\catcode`\<=\active \def<{
\fontencoding{T1}\selectfont\symbol{60}\fontencoding{\encodingdefault}}
\newcommand{\equallim}{\mathop{=}\limits}
\newcommand{\tmaffiliation}[1]{\\ #1}
\newcommand{\tmem}[1]{{\em #1\/}}
\newcommand{\tmemail}[1]{\\ \textit{Email:} \texttt{#1}}
\newcommand{\tmstrong}[1]{\textbf{#1}}
\newcommand{\tmverbatim}[1]{\text{{\ttfamily{#1}}}}
\newcommand{\tmarxiv}{\textbf{arXiv subject classification:} }
\newcommand{\tmkeywords}{\textbf{Keywords:} }

\begin{document}

\title{Argumentaci{\'o}n de futuros profesores de matem{\'a}ticas en tareas
sobre fracciones mediadas por un sistema de evaluaci{\'o}n en l{\'i}nea con
feedback autom{\'a}tico}

\author{
  Jorge Gaona
  \tmaffiliation{Universidad Metropolitana de Ciencias de la Educaci{\'o}n}
  \tmemail{jorge.gaona@umce.cl}
  \and
  Romina Menares
  \tmaffiliation{Universidad de Valpara{\'i}so}
  \tmemail{romina.menares@uv.cl}
}

\maketitle

\begin{abstract}
  {\tmstrong{Resumen{\tmem{{\tmstrong{\tmverbatim{}}}}}}}
  
  En este art{\'i}culo se reporta el trabajo de argumentaci{\'o}n de un grupo
  de profesores de matem{\'a}ticas en formaci{\'o}n en una experimentaci{\'o}n
  realizada en una clase virtual (a causa de la emergencia del COVID-19),
  durante el 2020. Se trabaj{\'o} con una tarea sobre fracciones en un sistema
  de evaluaci{\'o}n en l{\'i}nea con preguntas con par{\'a}metros aleatorios e
  infinitas respuestas correctas posibles. A partir de esto, se desarroll{\'o}
  un trabajo de discusi{\'o}n sobre las estrategias y las justificaciones,
  argumentaciones y validaciones de estas estrategias y de otras conjeturas
  que surgieron. En este art{\'i}culo se analiza este trabajo desde un enfoque
  cualitativo usando como marco te{\'o}rico el Espacio de Trabajo
  Matem{\'a}tico. Los resultados evidencian que el trabajo de discusi{\'o}n
  provoc{\'o} que los profesores en formaci{\'o}n pudieran hallar
  interpretaciones para los algoritmos procesados por la computadora,
  potenciando el discurso y la argumentaci{\'o}n epist{\'e}mica en contexto de
  uso de artefactos tecnol{\'o}gicos. A su vez, los mismos discursos
  permitieron a los futuros profesores instrumentalizar los procesos para ser
  utilizados en nuevas tareas.
\end{abstract}

\tmarxiv{math.HO}

\tmkeywords{tecnolog{\'i}a digital, justificaci{\'o}n y\.{}prueba,
formaci{\'o}n de profesores, clases en l{\'i}nea. }

\section{Introducci{\'o}n}

Seg{\'u}n Balacheff (1987), hay que diferenciar entre una explicaci{\'o}n o
argumentaci{\'o}n, una prueba, y una demostraci{\'o}n. Una argumentaci{\'o}n
es un discurso que busca hacer entendible o convencer sobre la veracidad de
una proposici{\'o}n; una prueba es una explicaci{\'o}n aceptada por una
comunidad; y una demostraci{\'o}n es un tipo de prueba con una estructura
espec{\'i}fica que se caracteriza por ser una cadena de enunciados organizados
con reglas determinadas. En el {\'a}mbito educativo se propone valorar
aquellas pruebas que puedan ayudar a explicar y justificar conjeturas
apoy{\'a}ndose en elementos matem{\'a}ticos (Hanna, 2001).

En Chile, las bases curriculares de secundaria proponen de forma expl{\'i}cita
la argumentaci{\'o}n como una de las habilidades que se deben desarrollar en
matem{\'a}ticas (MINEDUC, 2019).

La comunidad de investigadores en educaci{\'o}n matem{\'a}tica del pa{\'i}s
ha estado reportando distintos resultados en torno a la argumentaci{\'o}n en
la formaci{\'o}n inicial y continua de profesores de matem{\'a}ticas. En
geometr{\'i}a, Nagel et al. (2008) declaran que muchos profesores en su primer
a{\~n}o de formaci{\'o}n inicial son capaces de ordenar sus argumentos de modo
deductivo, sin embargo, a{\'u}n no poseen la capacidad de desarrollar
argumentos deductivos en todas las tareas propuestas, a{\'u}n cuando se trata
de una universidad que recibe a estudiantes de alto nivel acad{\'e}mico. En
otra investigaci{\'o}n en formaci{\'o}n inicial sobre sistemas de ecuaciones
lineales, Rodr{\'i}guez-Jara et al. (2019) observan que los estudiantes
argumentan err{\'o}neamente sobre la base de la proporcionalidad de los
coeficientes, y seg{\'u}n los autores, esto muestra que no han coordinado
algunos procesos relacionados con los objetos matem{\'a}ticos involucrados.

Hay otro grupo interesante de resultados sobre las argumentaciones, en este
caso, de profesores en ejercicio. En educaci{\'o}n primaria, Pizarro et al.
(2018) concluyen que las debilidades de profesores en el conocimiento de
estimaci{\'o}n no permiten argumentar de forma adecuada. En otro trabajo sobre
n{\'u}meros enteros, Solar (2018) muestra la importancia de algunas
estrategias para que el profesor fomente la argumentaci{\'o}n: dar
oportunidades de participaci{\'o}n, gestionar el error y realizar ciertos
tipos espec{\'i}ficos de preguntas, como preguntas que favorezcan la
explicaci{\'o}n, evitar preguntas ret{\'o}ricas, hacer contra preguntas y
preguntas que mantengan el foco en la discusi{\'o}n.

En cuanto a la educaci{\'o}n secundaria, varias investigaciones dan cuenta de
las argumentaciones en el trabajo matem{\'a}tico propuesto por profesores
debutantes en sus clases. Por ejemplo, en el trabajo de Henr{\'i}quez-Rivas y
Montoya-Delgadillo (2015) o de Montoya-Delgadillo el al. (2014) se se{\~n}ala
que, en las clases de profesores debutantes sobre {\'a}lgebra y geometr{\'i}a,
el trabajo de argumentaci{\'o}n es relegado o est{\'a} casi ausente.
Adem{\'a}s, en esta misma l{\'i}nea y considerando la argumentaci{\'o}n como
parte del trabajo discursivo dentro de una actividad matem{\'a}tica m{\'a}s
general, que considera aspectos semi{\'o}ticos e instrumentales (Kuzniak et
al., 2016), Henr{\'i}quez-Rivas y Montoya-Delgadillo (2016) muestran que hay
tareas que fomentan lo discursivo, aunque a{\'u}n se observan dificultades
para desarrollar esta dimensi{\'o}n cuando se implementan. En el {\'a}rea de
las probabilidades, Montoya-Delgadillo et al. (2016) reportan que los
profesores debutantes se encuentran con que no han tenido la experiencia
necesaria para promover procesos argumentativos, concluyendo que para los
profesores, los procesos de demostraci{\'o}n no son viables en el nivel
escolar del cual se ocupan, lo que los lleva a privilegiar la operatoria,
sumado adem{\'a}s a que, por encontrarse muy cercano a{\'u}n a su
instituci{\'o}n formadora, su desempe{\~n}o en el aula est{\'a} tensionado
entre lo que piensan como matem{\'a}ticos y lo que piensan como profesores.

A nivel internacional, al hacer una b{\'u}squeda sobre la argumentaci{\'o}n
de profesores en la base bibliogr{\'a}fica de Scopus y Web of Science,
encontramos que la argumentaci{\'o}n en la formaci{\'o}n de profesores
tambi{\'e}n ha sido abordada en diferentes lugares. En Estados Unidos y Reino
Unido se estudi{\'o} a profesores de primaria en trabajos sobre aritm{\'e}tica
y se concluye que los profesores tienen una buena comprensi{\'o}n de la
distinci{\'o}n entre prueba y argumentos emp{\'i}ricos (Stylianides \&
Stylianides, 2009), aunque los profesores se deben enfrentar a una serie de
desaf{\'i}os cuando resuelven una tarea que implica argumentar, tales como,
diferenciar entre describir y explicar, o coordinar las interpretaciones con
las estrategias (Lo et al., 2008). Cuando est{\'a}n en el papel de profesores,
las mayores dificultades est{\'a}n relacionadas con implementar tareas
relacionadas con el razonamiento y la demostraci{\'o}n y gestionar los
h{\'a}bitos mentales preexistentes de los estudiantes que no estaban en
sinton{\'i}a con la creaci{\'o}n de un sentido matem{\'a}tico (Stylianides et
al., 2013).

Belin y Akar (2020) abordaron la argumentaci{\'o}n con n{\'u}meros reales y
mostraron que luego de una instrucci{\'o}n donde se trabajaron distintas
representaciones, los futuros profesores generaron sus propios ejemplos y
dibujaron diagramas mientras explicaban los enunciados y las notaciones dadas
en las argumentaciones, donde los decimales se utilizaban para trabajar con
n{\'u}meros peri{\'o}dicos. Este tipo de instrucci{\'o}n puede superar
obst{\'a}culos epistemol{\'o}gicos (Bachelard, 1938) sobre la comparaci{\'o}n
entre n{\'u}meros mediante distintas notaciones peri{\'o}dicas (Mena-Lorca et
al., 2014).

Cuando en la b{\'u}squeda bibliogr{\'a}fica se agregan palabras relacionadas
con la tecnolog{\'i}a para potenciar la argumentaci{\'o}n en futuros
profesores, las investigaciones son m{\'a}s bien escasas; est{\'a}n los
art{\'i}culos que utilizan software de matem{\'a}ticas para asistir la
argumentaci{\'o}n y hay otros donde la tecnolog{\'i}a es el medio de
comunicaci{\'o}n para dar soporte al trabajo de argumentaci{\'o}n de los
estudiantes. Los trabajos de Stupel y Ben-Chaim (2017) y de Zengin (2017)
utilizaron Geogebra para asistir el trabajo de argumentaci{\'o}n y prueba en
tareas de geometr{\'i}a. En Stupel y Ben-Chaim (2017) se presentan varios
teoremas distintos cuya demostraci{\'o}n se realiza inicialmente mediante el
software de geometr{\'i}a din{\'a}mica y posteriormente de forma cl{\'a}sica.
En Zengin (2017), a partir de una tarea se les mostr{\'o} a los estudiantes
m{\'u}ltiples soluciones, algunas de ellas se exploraban inicialmente con
Geogebra. En ambas investigaciones, los estudiantes tuvieron una buena actitud
hacia la argumentaci{\'o}n luego hacia el trabajo propuesto por los
investigadores.

Por su parte, Fern{\'a}ndez et al. (2012), presenta un trabajo donde se
utiliza la tecnolog{\'i}a como medio de comunicaci{\'o}n. En este trabajo se
subraya que el formato ciertamente presenta desaf{\'i}os adicionales al
proceso de argumentaci{\'o}n. La investigaci{\'o}n muestra una tarea de
aritm{\'e}tica, y se reflexiona acerca de c{\'o}mo se puede llevar a cabo una
discusi{\'o}n en l{\'i}nea tomando en cuenta que compartir un texto escrito
con otros puede fomentar la discusi{\'o}n.

En Chile, como consecuencia de la pandemia, las clases se hicieron en formato
en l{\'i}nea durante todo el a{\~n}o 2020, y la participaci{\'o}n de los
estudiantes fue una de las problem{\'a}ticas que emergieron. El acceso y la
conectividad para la asistencia a las clases virtuales ha sido la primera
barrera. Se ha reportado un promedio de 27\% para el quintil m{\'a}s pobre y
89\% para el quintil m{\'a}s rico de acceso a las clases virtuales en
educaci{\'o}n obligatoria (MINEDUC, 2020).

En educaci{\'o}n superior a{\'u}n no hay datos reportados. Sin embargo, a
nivel local, se ha detectado poca participaci{\'o}n de los estudiantes en las
clases virtuales, poco compromiso de los estudiantes para participar, y
c{\'a}maras y micr{\'o}fonos apagados durante las clases.

Dado este cambio obligado por la emergencia, en este art{\'i}culo analizamos
las argumentaciones epist{\'e}micas que pueden existir entre los estudiantes y
el profesor en el contexto de una clase virtual apoyada por un sistema de
evaluaci{\'o}n en l{\'i}nea para una tarea de fracciones, donde la
tecnolog{\'i}a aparece como elemento de comunicaci{\'o}n, entre estudiantes y
profesores, y tambi{\'e}n, como soporte de la actividad matem{\'a}tica. De
forma m{\'a}s espec{\'i}fica, la pregunta de investigaci{\'o}n que se propone
es: {\tmstrong{{\textquestiondown}cu{\'a}l es el trabajo de argumentaci{\'o}n
puesto en juego en un contexto de clase virtual donde se usa una tarea en un
artefacto de evaluaci{\'o}n en l{\'i}nea para una tarea de fracciones?}}

\section{Marco te{\'o}rico}

Queremos caracterizar las discusiones epist{\'e}micas usando el Espacio de
Trabajo Matem{\'a}tico (ETM) (Kuzniak et al., 2016a; Kuzniak \& Richard, 2014)
descomponi{\'e}ndolas en las distintas g{\'e}nesis y planos que ofrece este
marco te{\'o}rico.

El ETM permite identificar y analizar el trabajo matem{\'a}tico de un sujeto,
tomando en cuenta aspectos epistemol{\'o}gicos y cognitivos del tema abordado
y de quien resuelve, respectivamente. Estos aspectos se articulan mediante las
g{\'e}nesis semi{\'o}tica, instrumental y discursiva (Figura 1). Ac{\'a} la
palabra g{\'e}nesis se utiliza en un sentido amplio y hace referencia, tanto
al comienzo de un proceso, como a su desarrollo e interacci{\'o}n entre los
polos del plano epistemol{\'o}gico y del plano cognitivo.

\begin{figure}[h!]
	\begin{center}
 \raisebox{0.0\height}{\includegraphics[width=0.8\textwidth]{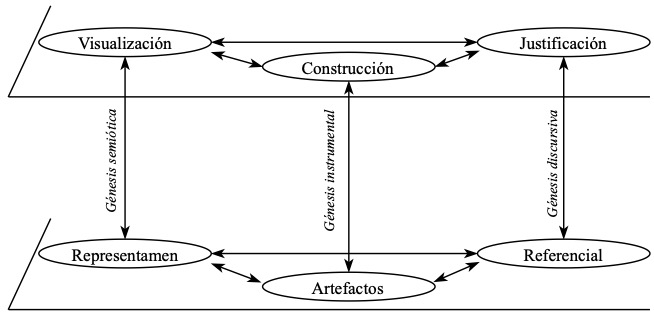}}
\caption{Espacio de trabajo matem{\'a}tico. Extra{\'i}do de Kuzniak (2016).}
	\end{center}
   
\end{figure}

En nuestro trabajo, para identificar el valor epist{\'e}mico de las
discusiones que se generan entre pares y estudiantes--profesor,
identificaremos cu{\'a}les son las g{\'e}nesis que se activan en las
argumentaciones que se realizan. Para esto, es indispensable identificar si un
objeto matem{\'a}tico, tal como una funci{\'o}n o una f{\'o}rmula, es usado en
el di{\'a}logo como una herramienta semi{\'o}tica, artefacto material o
simb{\'o}lico o como herramienta te{\'o}rica.

En el plano epistemol{\'o}gico se encuentran los objetos y/o herramientas que
permiten desarrollar el trabajo matem{\'a}tico y se definen tres polos: el
representamen , los artefactos y el referencial te{\'o}rico. De acuerdo con
Kuzniak et al. (2016b), en el modelo de los ETM, los objetos matem{\'a}ticos
pueden convertirse en herramientas o viceversa. Por otra parte, en el plano
cognitivo, se encuentran tres procesos, a trav{\'e}s de los cuales se intenta
dar cuenta de la actividad matem{\'a}tica: visualizaci{\'o}n, construcci{\'o}n
y prueba.

Cabe observar que, para poder identificar las g{\'e}nesis que se activan o
privilegian en la resoluci{\'o}n de una tarea, necesitamos identificar en
cu{\'a}l de los tres polos del plano epistemol{\'o}gico se encuentra un objeto
matem{\'a}tico en particular. El estatus de un objeto o herramienta, en
relaci{\'o}n con el polo en el cual se encuentra, estar{\'a} dado m{\'a}s bien
por su utilizaci{\'o}n que por una caracter{\'i}stica intr{\'i}nseca. A saber,
diremos que un objeto o herramienta matem{\'a}tica est{\'a} en el polo del
representamen cuando se utiliza como una herramienta semi{\'o}tica, en otros
t{\'e}rminos, cuando se trabaja a partir de su visualizaci{\'o}n y se toman en
cuenta las relaciones entre sus unidades figurales y no solo la percepci{\'o}n
visual que provee el acceso directo al objeto. Los objetos del polo de los
artefactos materiales o simb{\'o}licos se identificar{\'a}n cuando se trabaje
con herramientas materiales (como regla y comp{\'a}s), herramientas
inform{\'a}ticas (como una calculadora CAS) o artefactos simb{\'o}licos, como
un algoritmo. Para los dos primeros casos, dada su naturaleza, son
f{\'a}cilmente identificables, en cambio, el artefacto simb{\'o}lico lo
identificaremos cuando un objeto matem{\'a}tico o un algoritmo se utilice como
una herramienta para obtener un resultado y no se tomen en cuenta sus
propiedades, vale decir, cuando su uso no est{\'e} apoyado en el referencial
te{\'o}rico. As{\'i}, asociaremos el estatus de simb{\'o}lico a un uso
totalmente naturalizado y rutinario, en el que no se discuta ni cuestione su
validez ni justificaci{\'o}n.

Finalmente, en el polo del referencial te{\'o}rico est{\'a}n las propiedades,
teoremas y axiomas que dan sustento al discurso matem{\'a}tico. Este polo no
se debe pensar solo como una colecci{\'o}n de propiedades, porque al dar
soporte a las justificaciones deductivas, debe estar organizado de forma
coherente y bien adaptado a las tareas que se les pide a los estudiantes que
resuelvan (Kuzniak et al., 2016b).

Tal como lo indican Kuzniak et al. (2016a), la matem{\'a}tica es ante todo una
actividad humana y no solo una lista de signos y propiedades. Por esta
raz{\'o}n, este modelo considera un segundo nivel centrado en el sujeto,
consider{\'a}ndolo como un sujeto cognitivo cuyos procesos mentales est{\'a}n
en interacci{\'o}n con el plano epistemol{\'o}gico a trav{\'e}s de una
actividad matem{\'a}tica espec{\'i}fica. \ Los tres procesos que se consideran
en el plano cognitivo son: la visualizaci{\'o}n, la construcci{\'o}n y la
prueba. \ La visualizaci{\'o}n, est{\'a} relacionada con la interpretaci{\'o}n
de signos y la construcci{\'o}n interna de la representaci{\'o}n de los
objetos y sus relaciones. La construcci{\'o}n, est{\'a} relacionada con la
utilizaci{\'o}n de artefactos (materiales o simb{\'o}licos), junto con
esquemas de uso para producir elementos tangibles como escritos o dibujos y
tambi{\'e}n para la observaci{\'o}n, exploraci{\'o}n y experimentaci{\'o}n
mediada por un artefacto. Finalmente, la prueba, est{\'a} relacionada con el
proceso de justificaci{\'o}n mediante herramientas te{\'o}ricas y no solamente
una validaci{\'o}n emp{\'i}rica, la cual se podr{\'i}a entender m{\'a}s como
el proceso de construcci{\'o}n antes descrito.

Los planos epistemol{\'o}gico y cognitivo se articulan mediante tres
g{\'e}nesis: semi{\'o}tica, instrumental y discursiva.

La g{\'e}nesis semi{\'o}tica conecta el proceso de visualizaci{\'o}n en el
plano cognitivo con el representamen en el epistemol{\'o}gico. Esta
g{\'e}nesis puede partir por el signo en el representamen que es interpretado
por el sujeto mediante la visualizaci{\'o}n. Tambi{\'e}n puede partir por el
sujeto que codifica y produce un signo.

La g{\'e}nesis instrumental conecta el proceso de construcci{\'o}n en el
plano cognitivo con el polo de los artefactos. Cuando se trabaja con
herramientas materiales, inform{\'a}ticas o simb{\'o}licas, la g{\'e}nesis
involucra dos procesos: el de instrumentalizaci{\'o}n y el de
instrumentaci{\'o}n (Coutat \& Richard, 2011). El primero, comprende la
emergencia y evoluci{\'o}n de los esquemas de uso del artefacto y la
utilizaci{\'o}n de las posibilidades que ofrece el artefacto. El segundo parte
desde el sujeto y es relativo a la emergencia y evoluci{\'o}n de los esquemas
de uso y de las acciones instrumentadas, su constituci{\'o}n, funcionamiento,
coordinaci{\'o}n, combinaci{\'o}n, inclusi{\'o}n y asimilaci{\'o}n de
artefactos nuevos a esquemas ya constituidos. El trabajo matem{\'a}tico
podr{\'i}a ser considerado rutinario si es que no se conecta con la
validaci{\'o}n y justificaci{\'o}n de los artefactos.

Finalmente, la g{\'e}nesis discursiva conecta el proceso de prueba con el
polo del referencial te{\'o}rico en el plano epistemol{\'o}gico y est{\'a}
asociado al proceso de razonamiento deductivo mediante teoremas y propiedades.
En este {\'u}ltimo caso, el foco est{\'a} puesto en las propiedades y
teoremas, por lo que se est{\'a} pensando en razonamientos que van m{\'a}s
all{\'a} de los visuales o instrumentales, pero que pueden ser desencadenadas
por estos.

Cuando no es posible distinguir qu{\'e} g{\'e}nesis est{\'a} siendo
privilegiada, el ETM se puede caracterizar mediante la conexi{\'o}n de dos
g{\'e}nesis, considerando algunos de los tres planos verticales:
semi{\'o}tico-instrumental, semi{\'o}tico-discursivo o el
instrumental-discursivo (Coutat \& Richard, 2011).

\section{Toma de datos y metodolog{\'i}a}

La presente investigaci{\'o}n se desarrolla a trav{\'e}s de un enfoque
cualitativo, dada la naturaleza del fen{\'o}meno en estudio. De acuerdo con
Creswell y Poth (2009), la investigaci{\'o}n cualitativa es un proceso
interpretativo de indagaci{\'o}n basado en distintas tradiciones
metodol{\'o}gicas que examina un problema humano o social. Desde esta
perspectiva, y entendiendo que la problem{\'a}tica est{\'a} enmarcada en el
{\'a}mbito de la Did{\'a}ctica de la Matem{\'a}tica como ciencia, resulta
apropiada una indagaci{\'o}n que permita articular un marco te{\'o}rico con
procedimientos espec{\'i}ficos para el an{\'a}lisis de datos. El an{\'a}lisis
corresponde a un caso m{\'u}ltiple (Stake, 2007) de tipo opin{\'a}tico (Ruiz,
2003), pues se trata de un curso al cual se tiene acceso.

En la Universidad donde se realiz{\'o} el estudio, los alumnos de primer
a{\~n}o de la carrera Pedagog{\'i}a en Matem{\'a}ticas tienen una asignatura
denominada TICs para el aprendizaje de las matem{\'a}ticas. Uno de los
objetivos propuestos fue propiciar en el curso la discusi{\'o}n y la
argumentaci{\'o}n sobre distintos conceptos matem{\'a}ticos utilizando
herramientas tecnol{\'o}gicas. Esta universidad es p{\'u}blica en Chile y
tiene una larga tradici{\'o}n en la formaci{\'o}n de profesores de
matem{\'a}ticas.

Se propuso una serie de tareas en una plataforma Moodle (http://moodle.org)
que ten{\'i}a incorporado el plugin Wiris (http://www.wiris.com), a la cual
cada estudiante ten{\'i}a acceso durante la clase, pues bastaba con tener
conexi{\'o}n a internet. Estas herramientas permitieron crear tareas con
infinitas respuestas con el fin de propiciar la discusi{\'o}n. Las tareas se
implementaron durante el segundo semestre del 2020, fueron 15 profesores en
formaci{\'o}n (de primer a{\~n}o) quienes trabajaron en un primer momento en
forma individual respondiendo en la plataforma durante 10 minutos y luego vino
una fase de discusi{\'o}n. La clase, al ser virtual, fue grabada y se
transcribieron los pasajes donde se produce una discusi{\'o}n epist{\'e}mica.

En el cuestionario que se analiza en este trabajo cada estudiante ten{\'i}a un
espacio (virtual) individual donde deb{\'i}a dar tres respuestas distintas a
la misma tarea. A partir del trabajo de Gaona (2020), el artefacto-tarea en
una plataforma se evaluaci{\'o}n en l{\'i}nea se descompone en:
\begin{itemize}
  \item El enunciado: en el que a su vez se puede identificar el tipo de tarea
  y los objeto matem{\'a}tico Gaona et al. (2021): en este caso la tarea la
  tarea consisti{\'o} en ingresar una fracci{\'o}n de la forma a/b, con a y b
  enteros, que est{\'e} entre $\dfrac{1}{4}$ y $\dfrac{2}{4}$ junto con un
  trozo de la recta real donde se ubicaban estos valores.
  
  \item El sistema de entrada: un editor de ecuaciones para ingresar las
  respuestas. 
\end{itemize}
\begin{figure}[h!]
\begin{center}
  \raisebox{0.0\height}{\includegraphics[width=0.99\textwidth]{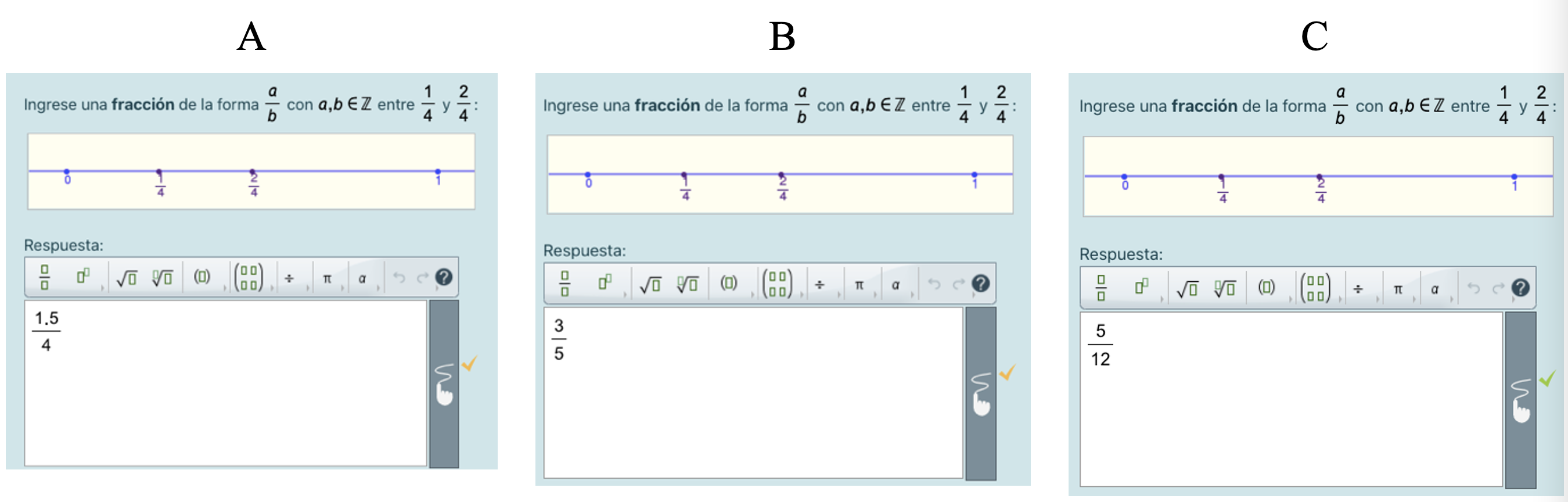}}
  \caption{Tarea que deb{\'i}a responder cada estudiante.}
\end{center}
\end{figure}
\begin{itemize}
  
  \item El sistema de validaci{\'o}n: que clasificaba las respuestas entre
  incorrectas, parcialmente correctas e incorrectas. El sistema califica la
  primera respuesta (figura 2A) como parcialmente correcta porque el
  estudiante ingres{\'o} un decimal en el numerador, y se ped{\'i}a en las
  instrucciones que fueran n{\'u}meros enteros. La segunda respuesta (figura
  2B) tambi{\'e}n fue calificada como parcialmente correcta porque el valor
  ingresado es mayor que $\dfrac{1}{4}$ pero no menor que $\dfrac{2}{4}$. La
  tercera respuesta (figura 2C) se considera correcta, pues cumple con todos
  los requisitos.
  
  \item El sistema de feedback: una vez que el estudiante ingresa la
  respuesta, el sistema adem{\'a}s de calificarlo como correcto o incorrecto,
  le entregaba un feedback autom{\'a}tico, que muestra por qu{\'e} las
  respuestas son consideradas como parcialmente correctas o incorrectas. Lo
  anterior se muestra en la figura 2.
  
  En este feedback, la plataforma ubica en la recta num{\'e}rica el valor
  ingresado y se indica si cumple con el resto de las condiciones, es decir,
  el feedback no da la respuesta ni sugiere alguna estrategia sobre c{\'o}mo
  resolverlo, sino que entrega una explicaci{\'o}n sobre por qu{\'e} lo
  realizado es correcto o incorrecto. Adem{\'a}s, en el caso de ser correcto,
  invita al usuario a buscar otro valor que sea distinto al ingresado, por
  ejemplo, en la figura 3C, se pide ingresar un valor m{\'a}s cercano a 1/4.
\end{itemize}

\begin{figure}[h!]
  \raisebox{0.0\height}{\includegraphics[width=0.99\textwidth]{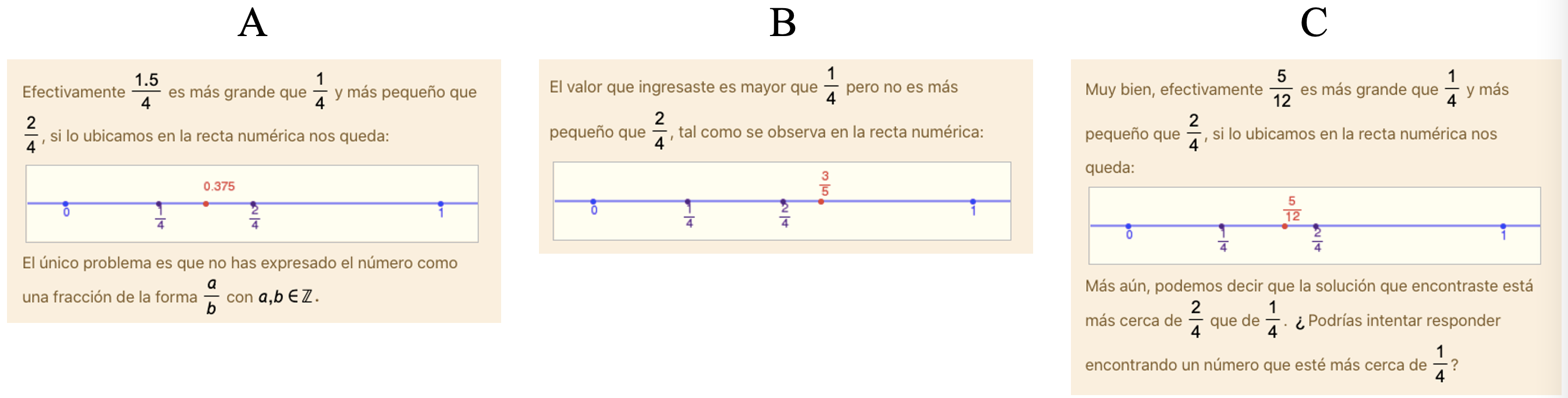}}
  \caption{Feedback autom{\'a}tico generado por el sistema una vez que el
  estudiante ingresa una respuesta.}
\end{figure}

En este caso, el artefacto digital contiene una componente pragm{\'a}tica y
otra epist{\'e}mica (Artigue, 2002): tiene una componente pragm{\'a}tica
porque recibe, registra y califica autom{\'a}ticamente la respuesta de los
estudiantes; y tiene una componente epist{\'e}mica porque emite un juicio
sobre la respuesta del estudiante (indicando si es correcto o incorrecto) y
porque, adem{\'a}s, en el feedback ubica la respuesta del estudiante en una
recta num{\'e}rica, haciendo comparaciones de forma gr{\'a}fica y analiza la
estructura de la respuesta para indicar si cumple o no con el formato
solicitado, por lo que da un significado a la interacci{\'o}n que se produce.

Durante el estudio, se trabaj{\'o} en tres etapas:
\begin{itemize}
  \item {\tmstrong{Clase 1:}} se hizo un trabajo individual en la plataforma
  de 10 minutos donde deb{\'i}an encontrar una fracci{\'o}n entre
  $\dfrac{1}{4}$ y $\dfrac{2}{4}$. El resto de la clase se discuti{\'o} sobre
  las estrategias y generalizaci{\'o}n de resultados. La clase termin{\'o} con
  la tarea: encuentre una secuencia infinita de fracciones entre
  $\dfrac{1}{4}$ y $\dfrac{2}{4}$.
  
  \item {\tmstrong{Clase 2:}} se hizo un trabajo grupal donde deb{\'i}an
  trabajar en buscar una respuesta a la tarea que qued{\'o} planteada al
  final.
  
  \item {\tmstrong{Evaluaci{\'o}n{\tmstrong{:}}}} se realiz{\'o} una parte
  virtual donde ten{\'i}an que encontrar fracciones entre otras dos fracciones
  dadas, de forma aleatoria y luego, en papel, \ ten{\'i}an que buscar una
  secuencia de infinitas fracciones que estuvieran entre las fracciones que
  hab{\'i}an respondido en la plataforma.
\end{itemize}
Para los an{\'a}lisis se codificaron los di{\'a}logos del profesor con la
letra P y la de los estudiantes con la letra E y un n{\'u}mero para
diferenciarlo: E1, E2, etc.

\section{Resultados y discusi{\'o}n}

\subsection{clase 1: encontrar una fracci{\'o}n entre dos fracciones}

\subsubsection{Trabajo individual en la plataforma}

Los 15 estudiantes ingresaron respuestas. El promedio de intentos, respuestas
correctas y tiempo utilizado se resume en la tabla 1.

	\begin{table}[h!]
		\begin{center}
		\begin{tabular}{lllll}
			\hline
			Intentos & \% R1 & \% R2 & \% R3 & Tiempo (min)\\
			\hline
			1.3 & 80.8\% & 91.6\% & 92\% & 6.2\\
			\hline
		\end{tabular}
		\caption{Promedio de intentos, tasa de respuesta 1, 2 y 3 y promedio de
			tiempo (en minutos)}
		\end{center}
	\end{table}

El promedio de intentos fue de 1,3 porque hubo 5 estudiantes que hicieron 2
intentos y el resto hizo 1. Tambi{\'e}n, se puede observar que la tasa de
respuestas correcta fue alta, es decir, las preguntas resultaron ser
suficientemente sencillas para los estudiantes. De hecho, la tasa m{\'a}s baja
del primer intento se dio m{\'a}s bien por dificultades de tipo instrumental
como por ejemplo, un estudiante que ingres{\'o} la respuesta y una frase, que
el sistema consider{\'o} incorrecta. El promedio de tiempo empleado,
tambi{\'e}n indica que pudieron responder a la pregunta de forma bastante
r{\'a}pida.

\subsubsection{Estrategias declaradas por los estudiantes}

Este trabajo, inicialmente sencillo, permiti{\'o}, en una segunda instancia,
plantear preguntas durante la clase que se direccionaban hacia la
elaboraci{\'o}n de argumentos sobre las distintas estrategias utilizadas. Al
preguntar, mediante el chat o mediante voz, fueron 2 estrategias que
aparecieron:
\begin{itemize}
  \item 7 estudiantes indicaron que utilizaron decimales, es decir,
  transformaron las fracciones a decimales, encontraron un decimal entre los
  dos valores y luego convirtieron en fracci{\'o}n nuevamente el n{\'u}mero
  encontrado.
  
  \item 3 estudiantes declararon que utilizaron amplificaci{\'o}n, es decir,
  amplificaron $\dfrac{1}{4}$ y $\dfrac{2}{4}$ por distintos n{\'u}meros
  enteros positivos y a partir de eso encontraron un valor entre ambos.
  
  \item 5 estudiantes no se pronunciaron sobre la estrategia utilizada.
\end{itemize}
Mediante un proceso inductivo, se les pregunt{\'o} cu{\'a}ntas fracciones se
obten{\'i}an entre $\dfrac{1}{4}$ y $\dfrac{2}{4}$ al amplificar por $2, 3, 4,
\ldots$ y $n$. Los estudiantes r{\'a}pidamente conjeturaron que hay $n - 1$
fracciones entre ambas. Se les pidi{\'o} demostrar esta conjetura y se propuso
revisarla m{\'a}s tarde.

Luego se pidi{\'o} a las estudiantes describir el procedimiento que hab{\'i}an
realizado con decimales, una estudiante por ejemplo eligi{\'o} 0.3, 0.26 y
0.4, los que transform{\'o} a fracciones: $\dfrac{3}{10}, \dfrac{13}{5}$ y
$\dfrac{2}{5}$ (los simplific{\'o}).

\subsubsection{Trabajo de discusi{\'o}n a partir de una estrategia no
prevista}

Posteriormente, se les pregunt{\'o} a los estudiantes si alguien hab{\'i}a
utilizado los decimales de una forma distinta y E1 indic{\'o} que los
hab{\'i}a escogido al azar; nombr{\'o}: $\dfrac{3}{10} , \dfrac{2}{5}$ y el
$\dfrac{9}{20}$. A partir de esto, se produce el di{\'a}logo que se muestra en
la Tabla 1.

En esta discusi{\'o}n se observa que el estudiante E1 utiliza como estrategia
el ensayo y error. La calculadora la utiliza como artefacto digital. En
t{\'e}rminos te{\'o}ricos, podemos evidenciar la activaci{\'o}n parcial de la
g{\'e}nesis instrumental (Kuzniak et al, 2016a). Decimos parcial, pues el
trabajo no inicia a partir de la instrumentalizaci{\'o}n del artefacto que se
pone en uso, sino m{\'a}s bien, es un artefacto que se asocia a la acci{\'o}n
de corroborar si los n{\'u}meros establecidos son v{\'a}lidos o no para
satisfacer el enunciado de la pregunta que se planteaba. En esta discusi{\'o}n
el estudiante est{\'a} intentando explicitar, a partir de las preguntas del
profesor cu{\'a}les fueron sus procedimientos.

\begin{table}[h!]
  \begin{tabular}{p{12.0cm}}
    \hline
    110. P: el $\frac{9}{20} ${\textquestiondown}c{\'o}mo lo encontraste?
    porque es distinto a los que estaban ac{\'a} [refiri{\'e}ndose a los que
    hab{\'i}a encontrado previamente su compa{\~n}era]\\
    111. E1: empec{\'e} a amplificar, o sea empec{\'e} a subir el numerador,
    por consecuente el denominador tambi{\'e}n lo fui subiendo hasta que me
    diera una fracci{\'o}n\\
    112. P: quiero que explicites m{\'a}s eso mismo que est{\'a}s diciendo, a
    qu{\'e} te refieres con{\ldots}\\
    113. E1: o sea, que empec{\'e} a probar con n{\'u}meros, empec{\'e} a
    jugar con los n{\'u}meros y{\ldots} fui probando con todos 1, 2, 3,... 7,
    8, 9 si igual como que us{\'e} harto rato y empec{\'e} a dividir por
    n{\'u}meros, a ver si me daban entre{\ldots}\\
    114. P: ya, entonces, por ejemplo, {\textquestiondown}probaste con otros
    denominadores? {\textquestiondown}c{\'o}mo decid{\'i}as si no te
    funcionaba el denominador que hab{\'i}as elegido?\\
    115. E1: porque lo ve{\'i}a con la calculadora (risas)\\
    116. P: dame un ejemplo, quiero que explicites el proceso que realizaste,
    antes del 20 te acuerdas qu{\'e} n{\'u}mero elegiste?\\
    117. E1: el 5 parece\\
    118. P: calculaste $\frac{9}{5}$, {\textquestiondown}eso?\\
    119. E1: s{\'i}, pero ese es 1,8 como que es mucho, a ver un denominador
    muy grande, por ejemplo el 10\\
    120. P: y ese {\textquestiondown}cu{\'a}nto te dio?\\
    121. E1: 0.9, tampoco cumpl{\'i}a, luego eleg{\'i} el 15 yendo como
    par{\'a}metro, porque si me sal{\'i}a muy cercano sab{\'i}a que pod{\'i}a
    probar entre esos dos n{\'u}meros nuevos [] luego el 20 [el denominador
    elegido] que es 0.45.\\
    \hline
  \end{tabular}
  \caption{Momento 1, se intenta explicitar una estrategia definida por el
  estudiante como ``al azar''.}
\end{table}

\begin{figure}[h!]
	\begin{center}
 \raisebox{0.0\height}{\includegraphics[width=0.8\textwidth]{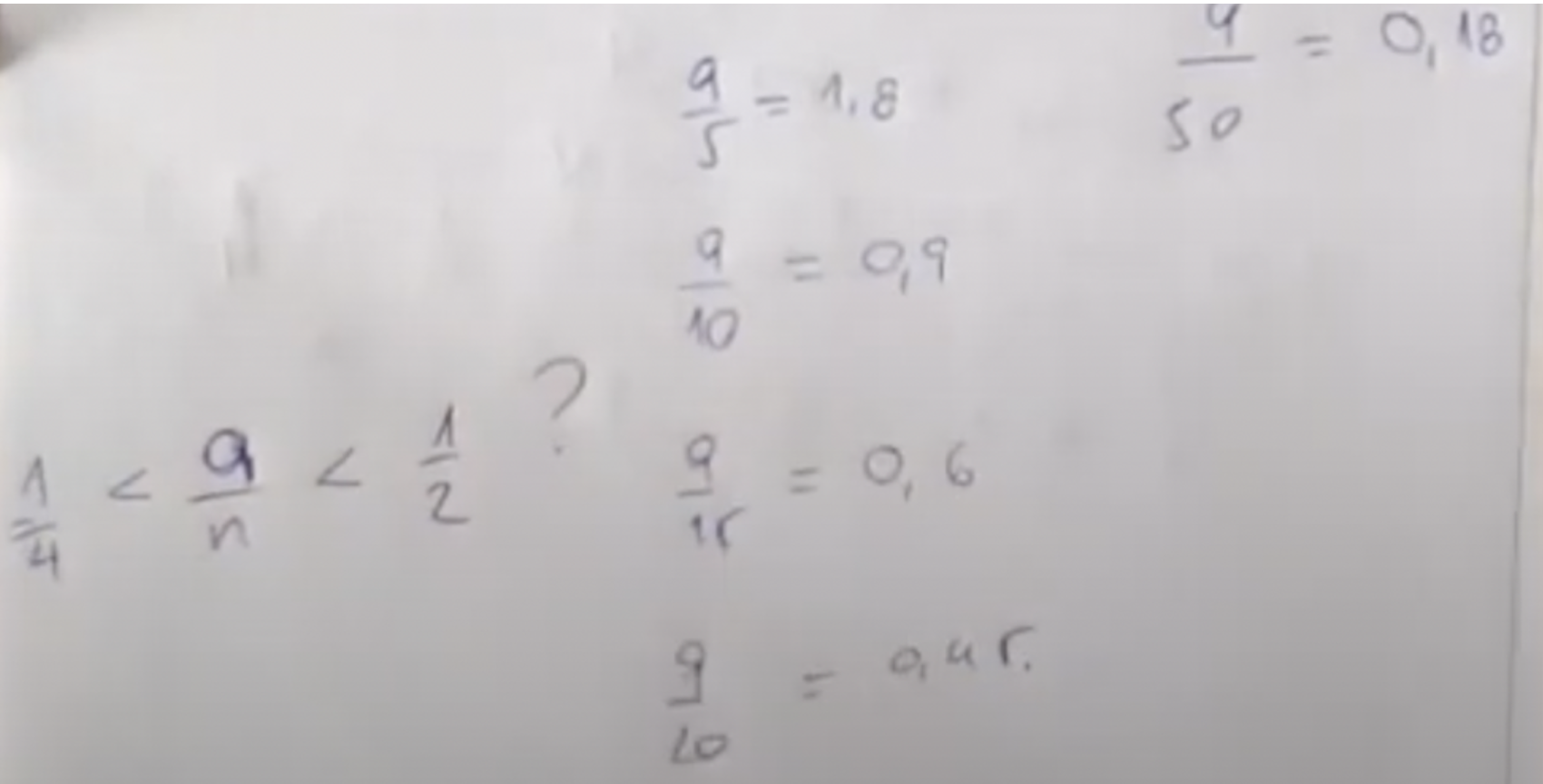}}
\caption{El profesor escribe lo que relata el estudiante E1 entre las
	l{\'i}neas 110 a 121, tambi{\'e}n escribe la pregunta que se indica en la
	l{\'i}nea 122 y 126 de la tabla 3}
	\end{center}
 
\end{figure}

El profesor, a partir del relato del estudiante sobre la estrategia de ensayo
y error, hace preguntas (l{\'i}nea 122 y 126 de la Tabla 2) que cambian el
foco de la discusi{\'o}n. Ahora se buscan los valores para los cuales
$\dfrac{9}{n}$ est{\'a}n entre $\dfrac{1}{4}$ y $\dfrac{2}{4}$. \ \

E2, en la l{\'i}nea 129, indica que ``elevando la fracci{\'o}n a menos 1 se
da vuelta la desigualdad''. \ En esta frase se observa que su argumento se
basa en un esquema rutinario, en el cual no explicita una propiedad, por lo
que calificamos su justificaci{\'o}n como instrumental. \

El profesor pregunta cu{\'a}les son las propiedades que justifican que se
inviertan las desigualdades. Es decir, nuevamente hay un cambio de foco. Si
primero se busc{\'o} explicitar las estrategias, luego encontrar todos los
valores de n que satisfac{\'i}an $\dfrac{1}{4} < \dfrac{9}{n} < \dfrac{2}{4}$,
ahora se busca, entre las l{\'i}neas 129 y 169, la justificaci{\'o}n de esta
soluci{\'o}n. Se observa que hay argumentos err{\'o}neos: como el de la
l{\'i}nea 133 donde E1 dice que se multiplica por $- 1$. O argumentos
circulares como el de la l{\'i}nea 140 dada por E4: ``se podr{\'i}a hacer
eliminando las fracciones // o sea pasando $a$ y $b$ a sus lados contrarios
(escrito por chat de la plataforma)'', el de la l{\'i}nea 141 dada por E1:
``{\textquestiondown}no puede multiplicar cruzado?'', E7 en la l{\'i}nea 150:
``se eleva a $- 1$ y se invierten las desigualdades?''.

Finalmente, entre las l{\'i}neas 155 y 169 las justificaciones se basan en el
referencial te{\'o}rico, los estudiantes evocan los axiomas y propiedades de
orden de los n{\'u}meros reales. Si uno compara los argumentos utilizados
antes de la l{\'i}nea 155 se observa un tr{\'a}nsito desde lo instrumental a
lo discursivo, los estudiantes pasan de usar artefactos simb{\'o}licos como
``multiplicar cruzado'' o ``se invierten'' a propiedades del referencial
te{\'o}rico, incluso, a diferencia de lo observado en el uso de artefactos
simb{\'o}licos, se preocupan del dominio de validez de la proposici{\'o}n
construida. Concretamente, se extiende la validez de la afirmaci{\'o}n m{\'a}s
all{\'a} de lo que se necesita para justifica la soluci{\'o}n de la
desigualdad.

En este proceso se observa que el profesor es quien valida los argumentos de
los estudiantes a trav{\'e}s de contra-preguntas. Hasta aqu{\'i}, los
estudiantes a{\'u}n no se muestran aut{\'o}nomos para determinar la validez de
sus afirmaciones.

Luego el profesor realiza un resumen donde retoma los elementos entregados por
los estudiantes y corrobora para $n = 21$: si $n$ es mayor a 18 y menor a 36
entonces $\dfrac{9}{21} = \dfrac{3}{7}$, que est{\'a} entre $\dfrac{1}{4}$ y
$\dfrac{2}{4}$. Primero, lo hace pidi{\'e}ndoles el c{\'a}lculo de su valor
decimal y luego les pide que demuestren que est{\'a} entre los valores
pedidos, sin usar decimales. En los c{\'a}lculos aparece el m{\'i}nimo
com{\'u}n m{\'u}ltiplo, que en una primera instancia obtuvieron por medio de
una tabla (artefacto simb{\'o}lico) y luego por descomposici{\'o}n prima. En
este episodio nuevamente se ve un tr{\'a}nsito de lo instrumental a lo
discursivo sobre el orden de las fracciones.

En seguida, se retoma una conjetura que se hab{\'i}a hecho al comienzo de la
clase y E2 la demuestra, lo que se transcribe en la Tabla 4. El estudiante
recurre a una secuencia escrita de forma algebraica para demostrar (Figura 4).
El argumento se transcribe en la l{\'i}nea 173. Al costado derecho, dentro de
la Figura 4, se ve tambi{\'e}n la demostraci{\'o}n discutida entre las
l{\'i}neas 122 y 169. Nuevamente se observa un trabajo en el plano discursivo
coordinado con artefactos simb{\'o}licos, donde el {\'a}lgebra sirve de
soporte para las demostraciones.

\begin{figure}[h!]
	\begin{center}
 \raisebox{0.0\height}{\includegraphics[width=0.8\textwidth]{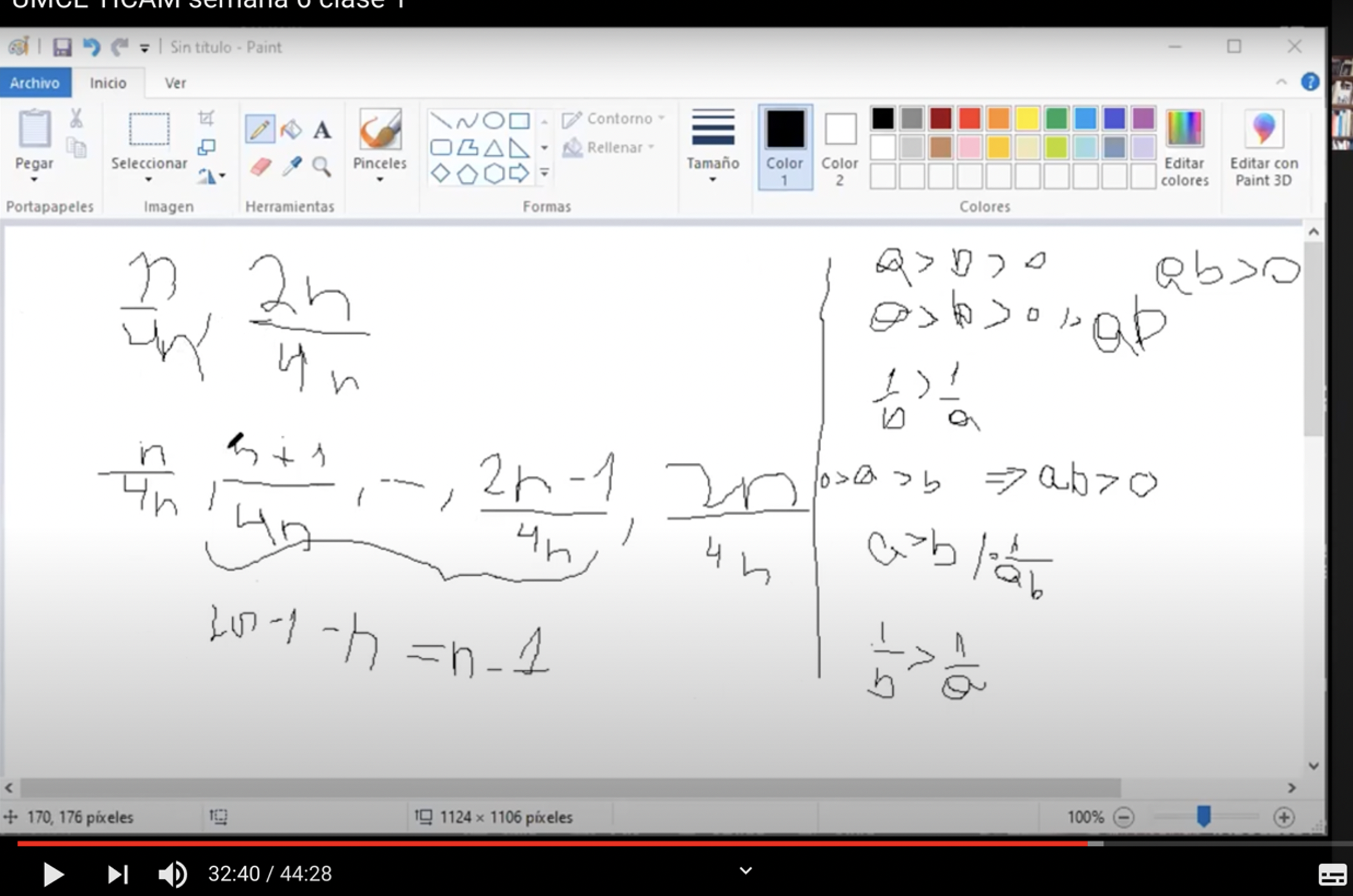}}
\caption{Pantalla compartida por E2 explicando la demostraci{\'o}n a la
	conjetura planteada.}
	\end{center}
 
\end{figure}

\newpage

\begin{table}[h!]
  \begin{tabular}{p{12.0cm}}
    \hline
    122. P: Si t{\'u} agrandas este n{\'u}mero tanto como t{\'u} quieras
    [refiri{\'e}ndose al denominador] {\textquestiondown}siempre est{\'a} en
    el intervalo?\\
    123. E1: creo que s{\'i}\\
    124. P: por ejemplo, piensa en un n{\'u}mero grande\\
    125. E1: $\frac{9}{50}$ da 0,18 no da dentro [del intervalo]\\
    126. P: el 9 hasta qu{\'e} n{\'u}mero lo puedes dividir?
    {\textquestiondown}$\frac{9}{n}$ es menor que $\frac{1}{2}$ y mayor que
    $\frac{1}{4}$? {\textquestiondown}alguien me podr{\'i}a dar la respuesta
    de eso? [el profesor va escribiendo parte de lo que relatan los
    estudiantes y de las preguntas que el hace ]\\
    127. E2: [despu{\'e}s de 20 segundos] entre 36 y 18\\
    128. P: {\textquestiondown}C{\'o}mo se obtiene ese 36 y 18?\\
    129. E2: puedo elevar toda la fracci{\'o}n a menos 1 y los s{\'i}mbolos se
    dan vuelta\\
    130. P: el profesor escribe eso en la hoja proyectada y pregunta
    {\textquestiondown}qu{\'e} pasa con los s{\'i}mbolos?\\
    131. E2: se invierten y amplifico todo por 9\\
    132. P: y queda $36 > n > 18$ {\textquestiondown}est{\'a}n de acuerdo?
    {\textquestiondown}por qu{\'e} esta desigualdad se invierte?\\
    {\tmstrong{133. E1: se multiplica por un n{\'u}mero nega{\ldots} osea, no
    se por qu{\'e}{\ldots}ahh sipo, se da vuelta porque multiplic{\'o} por $-
    1$}}\\
    134. P: no multiplic{\'o} por $- 1$, elev{\'o} a $- 1$\\
    135. E1: ahh ya, entonces no se por qu{\'e}\\
    136. P: la pregunta es, la voy a plantear de esta manera: Si $a < b$
    entonces {\textquestiondown}$\frac{1}{a} > \frac{1}{b}$?\\
    137. E3: s{\'i}, yo creo que s{\'i}\\
    138. P: entonces, esto es cierto para todo $a, b$ en $\mathbb{R}$ (los
    reales)? {\textquestiondown}C{\'o}mo lo podr{\'i}as demostrar? [...]
    {\textquestiondown}c{\'o}mo se demuestra? si es cierto lo que hizo su
    compa{\~n}ero antes (se los voy a mostrar) eso que est{\'a} ah{\'i} es
    cierto o sabr{\'i}amos en qu{\'e} casos lo podemos utilizar
    {\textquestiondown}ideas?\\
    139. E4: inaudible, [estudiantes indican que se escucha despacio e indica
    que escribir{\'a} por chat]\\
    {\tmstrong{\begin{tabular}{p{12.0cm}}
      140. E4 (por chat): se podr{\'i}a hacer eliminando las fracciones // o
      sea pasando $a$ y $b$ a sus lados contrarios\\
      141. E1: {\textquestiondown}no puede multiplicar cruzado?
    \end{tabular}}}\\
    142. E4. lo que dijeron reci{\'e}n jajajaja\\
    143. P: {\textquestiondown}se puede multiplicar cruzado? O sea pasando $a$
    y $b$ a sus lados contrario. T{\'u} dices (escribe $a < b$)\\
    144. E2: profe ya lo demostr{\'e}\\
    145. P: {\textquestiondown}cu{\'a}l, esta o la anterior?\\
    146. E2: las dos\\
    147. P: esp{\'e}rame un poquito para ver si el argumento que est{\'a}n
   utilizando o que est{\'a}n pensando los compa{\~n}eros es el mismo que vas
   a utilizar t{\'u}, si es distinto bien, si es la misma me dices si
   est{\'a}s de acuerdo o no. Si $a < b$ {\textquestiondown}qu{\'e}
   operaci{\'o}n debo aplicar a ambos lados de la desigualdad?\\
   148. E5: por $\frac{1}{a}$ y luego $\frac{1}{b}$\\
    \hline
  \end{tabular}

\end{table}

\begin{table}[h!]
	\begin{tabular}{p{12.0cm}}
		\hline
		
		149. E6: para eliminar denominadores\\
		{\tmstrong{150. E7: se eleva a -1 y se invierten las desigualdades?}}\\
		151. E2: eso es lo que queremos demostrar xd\\
		152. P: pero eso es lo que queremos demostrar, acu{\'e}rdense que as{\'i}
		parti{\'o} esto, elevo a menos 1 y se invierte la desigualdad\\
		153. E7: ah chuta no hab{\'i}a cachado\\
		154. P: est{\'a}n de acuerdo que demostrar ``$a < b$ entonces
		$\frac{1}{a}$<$\frac{1}{b}$'' es equivalente a demostrar ``si $a < b$ y
		elevo a la menos 1 entonces la desigualdad se invierte''?\\
		{\tmstrong{155. E1: {\textquestiondown}si igualamos a cero?
				{\textquestiondown}si $a < b$ entonces $a - b < 0$?}}\\
		156. P: entre este paso y el otro paso {\textquestiondown}qu{\'e}
		operaci{\'o}n hubo?\\
		{\tmstrong{157. E7: $- b$}}\\
		158. E1: se sum{\'o} el inverso aditivo\\
		159. P: {\textquestiondown}por qu{\'e} yo se que se mantiene la
		desigualdad? {\textquestiondown}Qu{\'e} axioma estoy utilizando?
		{\textquestiondown}Qu{\'e} axioma de orden me permite hacer este paso que
		hay ac{\'a}?\\
		{\tmstrong{160. E2: $a < b$ lo multiplico por $\frac{1}{ab}$ con $a$ y $b$
				positivos, as{\'i} me aseguro que no me invierta la desigualdad}}\\
		161. P: asumiste que $a$ y $b$ son positivos\\
		{\tmstrong{162. E2: s{\'i}, los restring{\'i} y queda
				$\frac{1}{b}$<$\frac{1}{a}$ es para los reales positivos}}\\
		163. P: hay est{\'a} lo que se quer{\'i}a demostrar\\
		{\tmstrong{164. E2: pero solo para reales positivos}}\\
		165. P: exactamente, el resultado no es para todo $\mathbb{R}$,
		{\textquestiondown}cierto? {\textquestiondown}hay otros n{\'u}meros para
		los cuales esto es cierto?\\
		166. E6: profe, {\textquestiondown}no ser{\'a} para cuando $a$ y $b$
		tienen el mismo signo? si $a$ y $b$ son negativos entonces tambi{\'e}n la
		multiplicaci{\'o}n da positiva y tienen que ser distinto de cero\\
		167. E1: si se puede en los reales positivos y negativos entonces
		{\textquestiondown}es en todos los reales? {\textquestiondown}o no?\\
		168. E6: deber{\'i}an ser en todos los reales siempre y cuando tengan el
		mismo signo\\
		169. P: entonces no es para cualquier par de reales, si alguno de los dos
		vale 0 {\textquestiondown}es cierto?\\
		\hline
	\end{tabular}
	\caption{Construcci{\'o}n de la justificaci{\'o}n de la soluci{\'o}n de la
		inecuaci{\'o}n $\frac{1}{4} < \frac{9}{n} < \frac{2}{4}$
	}
\end{table}

Para finalizar la clase el profesor plantea la siguiente pregunta:
{\textquestiondown}c{\'o}mo encontrar infinitas fracciones entre $\frac{1}{4}$
y $\frac{2}{4}$?

Esquematizamos lo ocurrido en la clase 1 en la Figura 6. En el esquema se
observa c{\'o}mo los estudiantes van transitando desde el uso de artefactos
materiales no cl{\'a}sicos (calculadora), a la utilizaci{\'o}n de artefactos
simb{\'o}licos (algoritmos). Luego, progresan en sus argumentos evocando al
referencial te{\'o}rico (con el uso de propiedades sobre n{\'u}meros reales).
As{\'i}, se construye una articulaci{\'o}n entre las g{\'e}nesis instrumental
y discursiva, es decir, una activaci{\'o}n del plano vertical
instrumental--discursivo, pues se aprecia una coordinaci{\'o}n entre procesos
que refieren al uso deartefactos, con argumentos basados en teor{\'i}a
matem{\'a}tica.

\begin{table}[h!]
  \begin{tabular}{p{12.0cm}}
    \hline
    170. P: E2 nos dijo que hab{\'i}a demostrado la otra conjetura. Que si
    tenemos $\frac{1}{4}$ y $\frac{2}{4}$ y lo amplifico por $n$ entonces hay
    $n - 1$ fracciones entre ambos n{\'u}meros y E2 dijo que lo hab{\'i}a
    demostrado. {\textquestiondown}Lo quieres mostrar?\\
    171. E2: profe lo hice en paint\\
    172. [...] P: expl{\'i}canos lo que estamos viendo (E2 est{\'a}
    demostrando que al amplificar por n se obtienen n-1 fracciones entre
    ambas)\\
    173. E2: defin{\'i} cuantas fracciones{\ldots} espere, me enred{\'e} Puse
    las fracciones una por una entre $\frac{n}{4 n}$ y $\frac{2 n}{4 n}$:
    $\frac{n}{4 n}, \frac{n + 1}{4 n}, \ldots ., \frac{2 n - 1}{4 n}, \frac{2
    n}{4 n}$ y cont{\'e} todas las fracciones que hab{\'i}an. Para contarlas
    rest{\'e} el numerador de ac{\'a} [refiri{\'e}ndose $2 n - 1$] y el de
    ac{\'a} [ refiri{\'e}ndose a $n$] y se obtiene $2 n - 1 - n$ y me dio $n -
    1$.\\
    174. P: Est{\'a} muy buena la demostraci{\'o}n, {\textquestiondown}alguien
    tiene alguna observaci{\'o}n?\\
    175. E2: puse todas las fracciones que hay entre $\frac{n}{4 n}$ y
    $\frac{2 n}{4 n}$ y luego las cont{\'e} calculando $2 n - 1 - n = n -
    1$.\\
    \hline
  \end{tabular}
  \caption{Transcripci{\'o}n de la demostraci{\'o}n de la conjetura: si
  amplifico el numerador y denominador de $\frac{1}{4}$ y $\frac{2}{4}$ por
  $n$, obtengo $n - 1$ fracciones entre ellas.}
\end{table}

\

\begin{figure}[h!]
  \raisebox{0.0\height}{\includegraphics[width=0.9\textwidth]{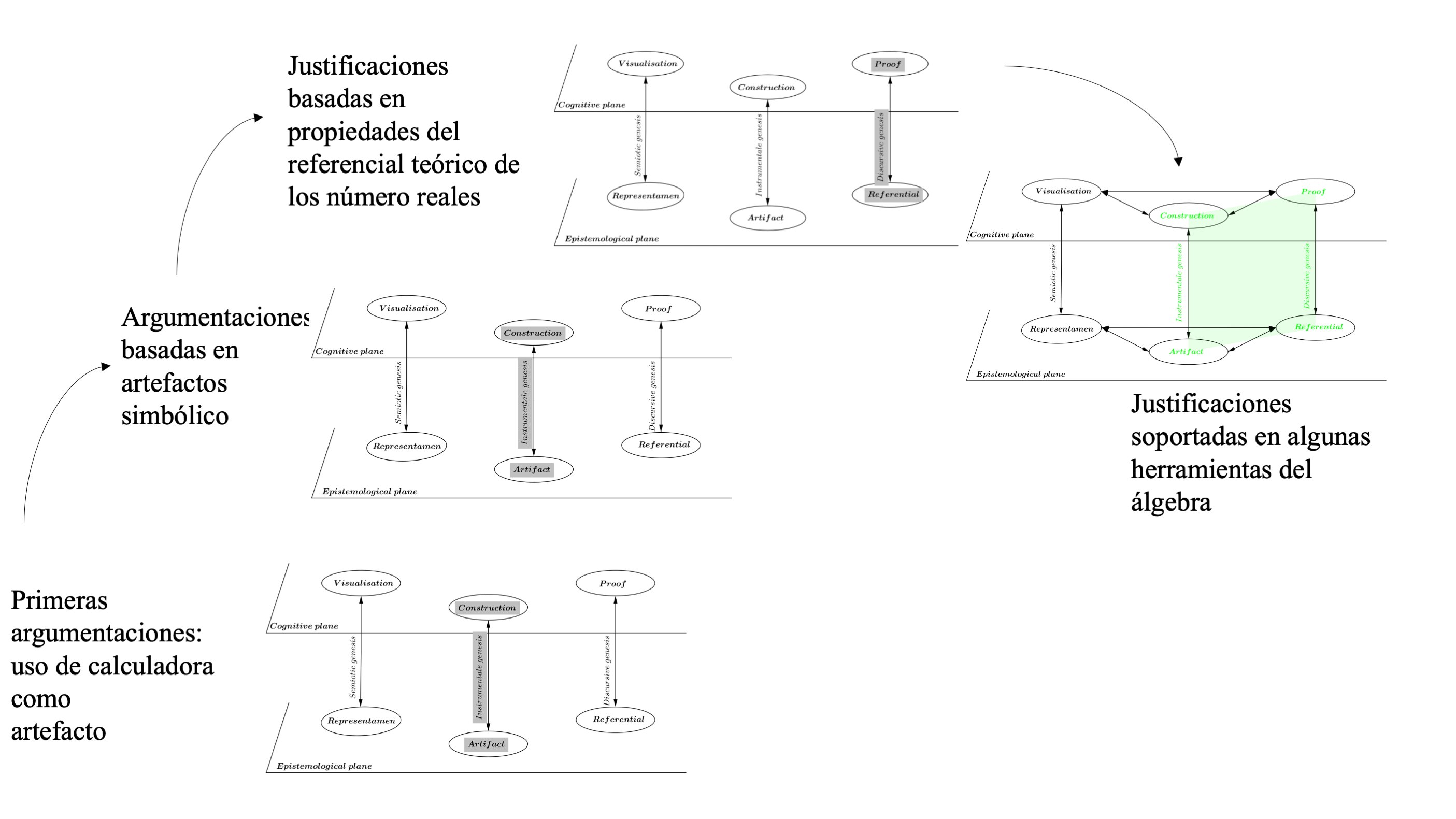}}
  \caption{Circulaci{\'o}n clase 1 a partir de la discusi{\'o}n sobre las
  estrategias utilizadas por los estudiantes para encontrar una fracci{\'o}n
  entre $\frac{1}{4}$ y $\frac{2}{4}$.}
\end{figure}

\subsection{Clase 2: infinitas fracciones entre $\frac{1}{4}$ y $\frac{2}{4}$}

En una segunda clase se abord{\'o} la pregunta planteada en la clase anterior,
para esto el curso se dividi{\'o} en 3 grupos, en el cual participaron 16
estudiantes.

\subsubsection{Trabajo en grupos}

\

El estudiante E1 del primer grupo explica en la l{\'i}nea 178 c{\'o}mo
construy{\'o} la secuencia (Tabla 5). Primero encontr{\'o} una secuencia
acotada entre 1 y 2, para luego dividirla por 4 para que quedara en el
intervalo solicitado. Ac{\'a} se utiliza de forma impl{\'i}cita que si $a <
a_n < b$ \ y $k > 0$ entonces $ka < ka_n < kb$. Esta argumentaci{\'o}n la
podemos clasificar en el plano instrumental-discursivo puesto que se
justifica, pero no est{\'a} claro hasta que punto el estudiante es consciente
de las propiedades que est{\'a} utilizando en su argumentaci{\'o}n. Se observa
que la secuencia que se obtiene es decreciente y converge a {\textonequarter},
pero esto no es mencionado. El profesor, realiza con los estudiantes un
trabajo de comprobaci{\'o}n con distintos n{\'u}meros donde la secuencia se
cumple. Los otros estudiantes del grupo muestran que no tienen claro que la
condici{\'o}n de que la secuencia cumpla la condici{\'o}n para ciertos valores
no es suficiente para que se cumpla para todo n mayor que cierto valor.

\begin{table}[h!]
  \begin{tabular}{p{12.0cm}}
    \hline
    178. E1: lo que entend{\'i} es que hab{\'i}a que encontrar una f{\'o}rmula
    general de n{\'u}meros entre $\frac{1}{4}$ y $\frac{2}{4}$. Lo primero que
    hice fue encontrar una que est{\'e} entre 2 y 1 y esta la divid{\'i} en 4
    que, ser{\'i}a $\frac{n + 1}{n}$ , luego lo divid{\'i} en 4 y qued{\'o}
    $\frac{n - 1}{4 n}$ est{\'a} entre $\frac{1}{4}$ y $\frac{2}{4}$.\\
    179. E2: entend{\'i} la primera parte, pero no entend{\'i} por qu{\'e} se
    divide por 4.\\
    180. P: primero veamos si es cierto lo que dice su compa{\~n}ero: $\frac{n
    - 1}{4 n}$ est{\'a} entre $\frac{1}{4}$ y $\frac{2}{4}$? Quiero que
    respondan E2 y E3.\\
    181. E3: si reemplazamos por 1 tendr{\'i}a que estar entre los valores\\
    182. P: tu dices si $n \equallim 1$ {\textquestiondown}cu{\'a}nto queda la
    expresi{\'o}n del medio?\\
    183. E3: $\frac{2}{4}$ que es lo mismo que $\frac{1}{2}$ y no cumple\\
    184. P: probemos otros n{\'u}meros porque lo que importa es que haya
    infinitos que sirvan. Si 4n={\guillemotright}
    {\textquestiondown}cu{\'a}nto queda?\\
    185. E3: 3/8\\
    186. P: {\textquestiondown}sirve?\\
    187. E3: s{\'i}, se comprueba multiplicando cruzado.\\
    188. P: probemos un n{\'u}mero m{\'a}s adelante: $n = 7$,
    {\textquestiondown}cu{\'a}nto queda?\\
    189. E3: $\frac{8}{28}$\\
    190. P: {\textquestiondown}ser{\'a} menor a $\frac{1}{2}$ y mayor a
    $\frac{2}{4}$?\\
    191. E3: se cumple\\
    192. P: {\textquestiondown}estos tres casos me sirven para demostrar que
    esta expresi{\'o}n est{\'a} siempre entre $\frac{1}{4}$ y $\frac{2}{4}$?\\
    193. E2: es que son poquitos n{\'u}meros\\
    194. P: si pruebo 10 y resulta {\textquestiondown}con eso ser{\'i}a
    suficiente?\\
    195. E3: multiplicando el n{\'u}mero que usemos en n nos va a dar un
    n{\'u}mero que sea menor que est{\'e} porque se va como achicando.\\
    196. P: si probamos 100, 1000 o hasta un mill{\'o}n no es suficiente para
    probar que todos est{\'a}n, necesitamos ese tipo de argumentos, un poco
    m{\'a}s espec{\'i}ficos para demostrar que todos est{\'a}n. Piensen en
    eso, me unir{\'e} a otro grupo que est{\'a} llamando.\\
    \hline
  \end{tabular}
  \caption{Transcripci{\'o}n del trabajo realizado por el grupo 1. Un
  estudiante encontr{\'o} una secuencia y sus compa{\~n}eros intentaban
  entender si funcionaba y por qu{\'e}.}
\end{table}

El grupo 2 (ver Tabla 6) manifiesta que no sabe c{\'o}mo abordar el problema.
En cambio, el grupo 3 indica que s{\'i} encontraron una, mostraron que
probaron que funcionaba porque evaluaron en n{\'u}meros grandes y
peque{\~n}os, de alguna forma utilizan lo que Balacheff (2000) denomina
``ejemplo crucial'' (p.26) porque buscan probar ejemplos generalizados. En el
di{\'a}logo se evidencia que la justificaci{\'o}n presentada se produce a
trav{\'e}s de la activaci{\'o}n del plano vertical instrumental-discursivo,
donde el proceso instrumental se da a trav{\'e}s de la evaluaci{\'o}n de
n{\'u}meros grandes en la secuencia construida. El profesor les indica que a
pesar de que, a partir de estos ensayos, pueden tener cierta seguridad de sus
resultados, esto no implica que este procedimiento pueda ser considerado como
una demostraci{\'o}n matem{\'a}tica.

\begin{table}[h!]
  \begin{tabular}{p{12.0cm}}
    \hline
    197. E4: profe no entendemos nada, estamos perdidos.\\
    198. E5: si profe, no sabemos c{\'o}mo hacerlo\\
    199. P: d{\'e}jenme ve a otros alumnos que est{\'a}n fuera de las grupos
    peque{\~n}os y vuelvo.\\
    200. P: {\textquestiondown}no saben c{\'o}mo hacerlo o no entienden el
    problema?\\
    201. E5: entiendo el problema, pero no se c{\'o}mo hacerlo, pens{\'e} en
    sumatorias, pero no me acuerdo.\\
    202. P: es m{\'a}s sencillo que eso, les dar{\'e} una idea para que puedan
    avanzar: piensen en $\frac{1}{n + 2}$ . Si $n = 1$, $\frac{1}{3}$ est{\'a}
    entre $\frac{1}{4}$ y $\frac{2}{4}$?\\
    203. E6: s{\'i}\\
    204. P: si $n = 2$, ${\textquestiondown} \frac{1}{4}$ est{\'a} entre
    $\frac{1}{4}$ y $\frac{2}{4} ?$\\
    205. E4, E5 y E6: no\\
    \hline
  \end{tabular}
  \caption{Transcripci{\'o}n de discusi{\'o}n del profesor con el grupo 2.}
\end{table}

\begin{table}[h!]
  \begin{tabular}{p{12.0cm}}
    \hline
    207. E7: encontramos una, pero restringido eso s{\'i}.\\
    208. P: {\textquestiondown}c{\'o}mo restringida?\\
    209. E7: solo para los n{\'u}meros naturales\\
    210. P: pero sirve, los n{\'u}meros naturales son infinitos. Habr{\'i}a
    que ver c{\'o}mo demostraron que estaba entre 1/4 y 2/4.\\
    211. E8: probamos en muchos n{\'u}meros\\
    212. E7: pero eso no es demostrar es mostrar\\
    213. E8: probamos con n{\'u}meros altos, con el mill{\'o}n, con el 100,
    probamos con n{\'u}meros grandes y chicos.\\
    214. P: {\textquestiondown}los probaron todos?\\
    215. E7: {\textquestiondown}c{\'o}mo los vamos a probar todos?\\
    216. P: no me convence su argumento, como probaron n{\'u}meros grandes y
    peque{\~n}os, eso les da cierta seguridad de que parece que resulta, pero
    como no los han probado todos [{\ldots}] eso es lo que tienen que
    demostrar\\
    217. E7: ya profe, ah{\'i} lo llamamos cuando est{\'e} listo.\\
    \hline
  \end{tabular}
  \caption{Transcripci{\'o}n de discusi{\'o}n del profesor con el grupo 3.}
\end{table}

\subsubsection{Discusi{\'o}n colectiva}

Luego hay una discusi{\'o}n colectiva con las secuencias encontradas. Solo
aparecen 3 secuencias, dos del grupo 1 y una el grupo 3, todas descritas en la
Tabla 8. Las que encontr{\'o} el grupo 1 son $a_n = \dfrac{2 n - 1}{4 n}$ y
$b_n = \dfrac{n + 1}{4 n}$, a pesar de que ambas secuencias sirven, E1 solo
pudo demostrar que a\_n no era mayor que {\textonequarter}. Este estudiante
utiliz{\'o} propiedades de los reales para ir acotando la secuencia y
demostrar lo que hab{\'i}a hecho, tal como se muestra la Figura 5.
Tambi{\'e}n, E7 del grupo 3 demuestra que la secuencia que encontr{\'o} es
creciente (ver Figura 6) aplicando la definici{\'o}n.

En ambas demostraciones se observa, a diferencia de lo que se hizo al comienzo
de la clase 1, que los estudiantes buscan cu{\'a}les son las propiedades y
definiciones que justifican sus proposiciones, es decir, este trabajo est{\'a}
m{\'a}s orientado hacia la g{\'e}nesis discursiva.

\begin{table}[h!]
  \begin{tabular}{p{12.0cm}}
    \hline
    218. P: E1 {\textquestiondown}puedes mostrar las secuencias que
    encontraste?\\
    219. E1: no las tengo juntas. Est{\'a} la otra secuencia $\frac{2 n - 1}{4
    n}$ y la demostraci{\'o}n para mayor o igual que $\frac{1}{4}$ me
    sali{\'o} mal.\\
    220. P: mu{\'e}strala igual porque, aunque no est{\'e} bien, puede dar
    pistas a tus compa{\~n}eros para que la intenten demostrar, porque lo que
    se demostr{\'o} es que $\frac{2 n - 1}{4 n}$ es mayor que $\frac{1}{4}$ y
    si se multiplica todo por $\frac{1}{n}$ no queda mayor que $\frac{1}{4}$.
    {\textquestiondown}y la otra que encontraste?\\
    221. E1: La otra, la demostraci{\'o}n es esta (muestra la Figura 5(c)) y
    la otra la tiene usted (el profesor hace una transcripci{\'o}n de lo
    relatado por E1 y lo muestra en la Figura 5(c))\\
    222. P: Lo interesante de esta secuencia es que tiene infinitos valores
    que se encuentran entre $\frac{1}{4}$ y $\frac{2}{4}$. Otro grupo propuso
    la secuencia $\frac{2 n - 1}{4 n}$, sobre esa secuencia hay varias
    preguntas: {\textquestiondown}est{\'a} entre $\frac{1}{4}$ y
    $\frac{2}{4}$? {\textquestiondown}es creciente, porque eso dijo el
    compa{\~n}ero?\\
    223. E7: demostr{\'e} que la segunda secuencia que hab{\'i}a puesto era
    creciente\\
    224. P: {\textquestiondown}puedes mostrar la demostraci{\'o}n?\\
    225. E7: no puse en escrito lo que hice pero lo que hice fue usar la
    definici{\'o}n que la diferencia entre $a_{n + 1} - a_n$ tiene que ser
    mayor o menor que cero. Si es mayor es creciente, si es menor es
    decreciente, entonces me dio un n{\'u}mero positivo as{\'i} que es
    creciente (muestra la Figura 6).\\
    226. P: dejamos hasta ac{\'a} la clase, estamos un poco pasados de tiempo,
    nos vemos la pr{\'o}xima semana.\\
    \hline
  \end{tabular}
  \caption{Transcripci{\'o}n de discusi{\'o}n colectiva, durante la clase 2,
  sobre las secuencias infinitas que se encontraron despu{\'e}s de trabajar en
  grupos}
\end{table}

\begin{figure}[h]
  \raisebox{0.0\height}{\includegraphics[width=0.99\textwidth]{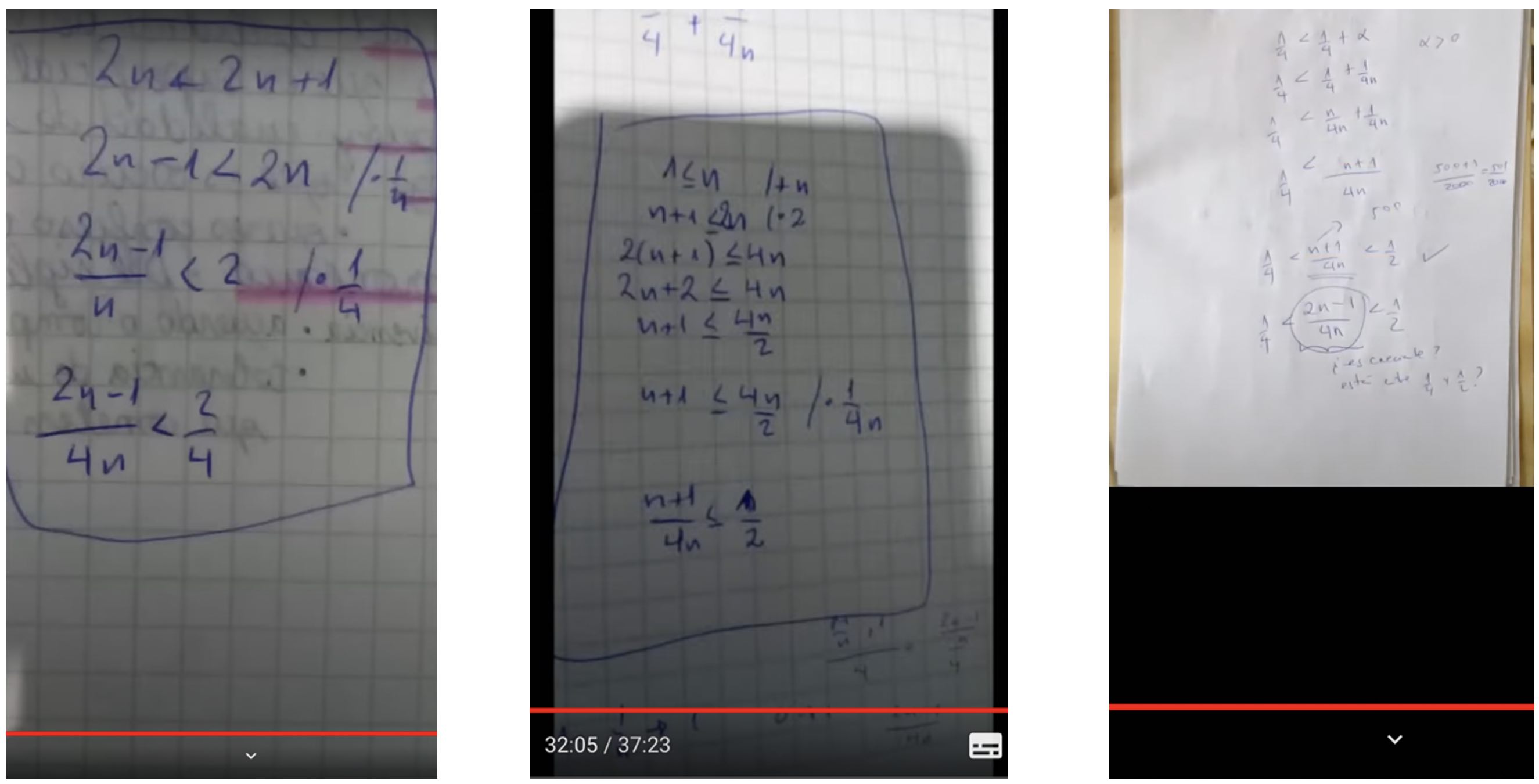}}
  \caption{Demostraci{\'o}n realizada por E1 en la clase 2 para justificar que
  las secuencias encontradas est{\'a}n acotadas entre $\frac{1}{4}$ y
  $\frac{2}{4}$.}
\end{figure}

\begin{figure}[h!]
	\begin{center}
	\raisebox{0.0\height}{\includegraphics[width=0.3\textwidth]{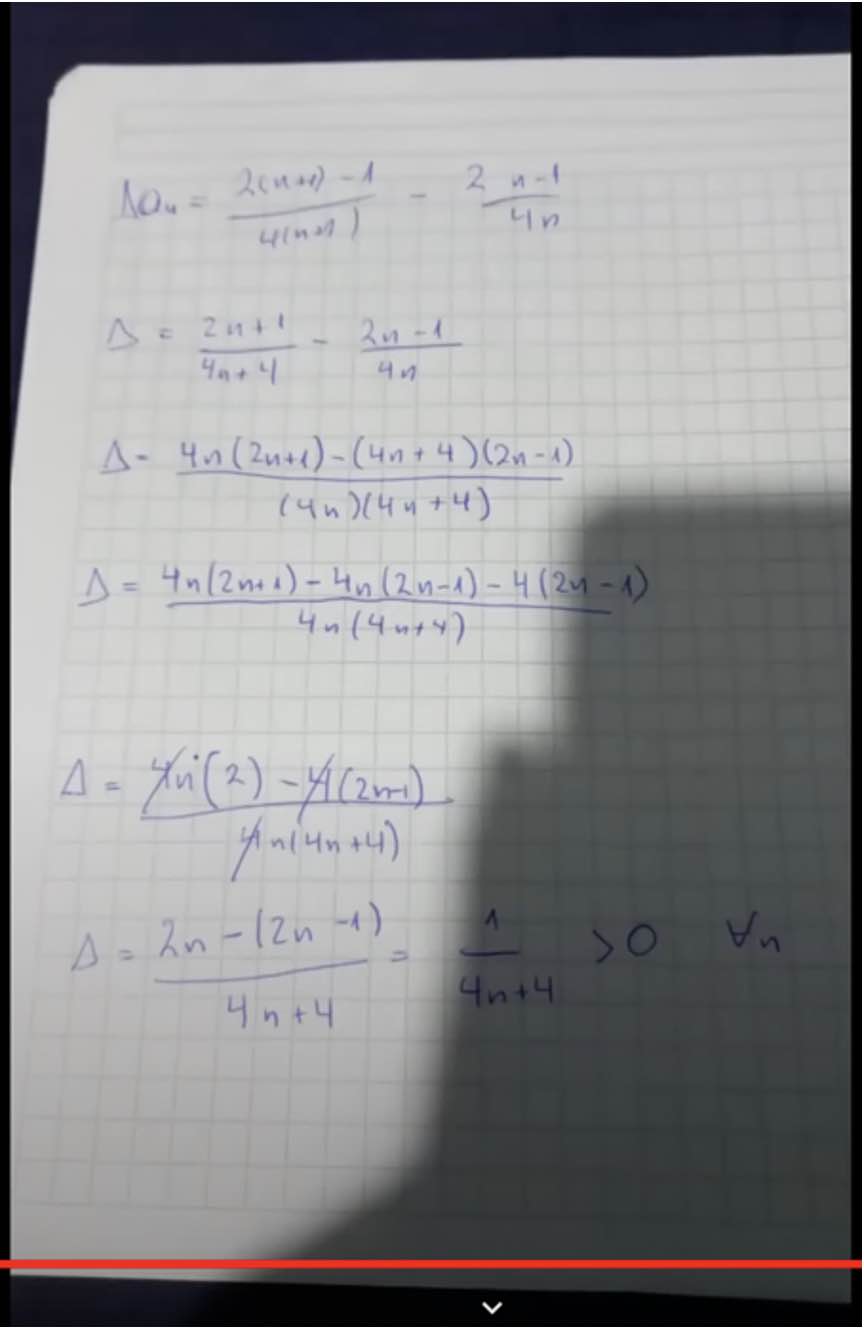}}
	\caption{Demostraci{\'o}n realizada por E1 en la clase 2 para justificar que
		una secuencia es creciente}
	\end{center}
\end{figure}

\

\

\

\

\subsection{Evaluaci{\'o}n de la secuencia}

Finalmente, se hizo una evaluaci{\'o}n que tuvo una componente individual en
la plataforma con correcci{\'o}n autom{\'a}tica y una grupal que deb{\'i}an
enviar en un escrito a mano o escrita en LaTex que se correg{\'i}a
manualmente. En esta evaluaci{\'o}n participaron 17 estudiantes, dos m{\'a}s
de los que asistieron a la primera clase. Todos dieron m{\'a}s de un ejemplo,
por eso hay m{\'a}s respuestas que n{\'u}mero de estudiantes.

Al hacer un resumen de todas las secuencias encontradas por los estudiantes,
en la Figura 7, se observa un patr{\'o}n: las respuestas que est{\'a}n entre
A1 y D5 (32 respuestas) construyen a partir de $\dfrac{a_1}{b_1}$ y
$\dfrac{a_2}{b_2}$, la secuencia $\dfrac{a_1}{b_1} + \dfrac{1}{b_1 n}$. Todas
las secuencias que se construyen son decrecientes y convergentes a
$\dfrac{a_1}{b_1}$.

El segundo grupo de respuestas corresponden a las tres secuencias que
est{\'a}n entre D6 y D8. En este caso, los estudiantes eligen una letra
distinta para la variable independiente. Construyen tres secuencias. La
primera, es una secuencia creciente y convergente al valor m{\'a}s grande del
intervalo, en este caso $\dfrac{6}{8}$. La segunda, es una secuencia como la
del primer grupo, es decir, creciente y convergente al valor m{\'a}s
peque{\~n}o, en este caso $\dfrac{1}{6}$. La tercera y {\'u}ltima secuencia
tambi{\'e}n es creciente, pero es una secuencia que est{\'a} conformada por
una potencia del valor m{\'a}s peque{\~n}o.

El tercer grupo son las secuencias de la columna E, que son todas secuencias
convergentes al centro del intervalo.

\begin{figure}[h]
\begin{center}
  \raisebox{0.0\height}{\includegraphics[width=0.9\textwidth]{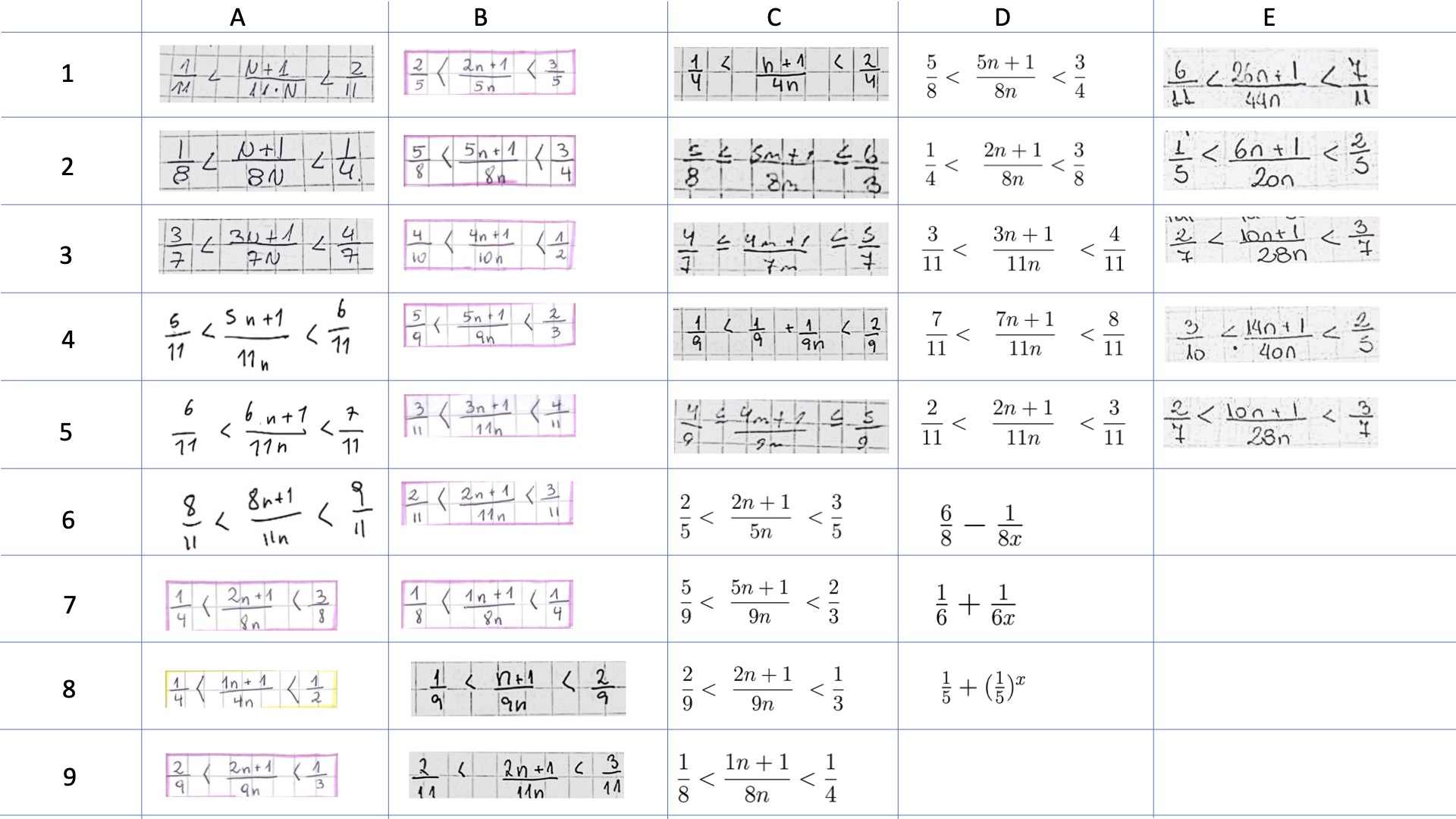}}
  \caption{Resumen de todas las secuencias encontradas por los estudiantes.}
\end{center}
\end{figure}

\subsubsection{Justificaciones del primer y tercer grupo de secuencias}

\

Las justificaciones del grupo 1 son demostraciones que siguen el mismo
esquema. Primero demuestran que $\dfrac{a_1}{b_1}$ es menor que
$\dfrac{a_1}{b_1} + \dfrac{1}{b_1 \cdot n}$, luego dejan la secuencia de la forma \
$\dfrac{a_1 \cdot n + 1}{b_1\cdot n}$ que es la respuesta que se presenta. Luego en una
segunda parte construyen, mediante desigualdades equivalentes, la cota
superior. Si se analiza a nivel individual (ver figura 10), la demostraci{\'o}n la
podr{\'i}amos clasificar en el plano discursivo, pero si se analiza a nivel
global, tomando en cuenta que la mayor{\'i}a de los estudiantes utilizaron el
mismo argumento adapt{\'a}ndolo a su secuencia, podr{\'i}amos concluir que
estas argumentaciones est{\'a}n en el plano instrumental-discursivo. De alguna
forma, el procedimiento para demostrar se instrumentaliza y pasa a convertirse
en un artefacto simb{\'o}lico, probablemente mediante un proceso de
discusi{\'o}n colectiva.

\

\begin{figure}[h!]
	\begin{center}
  \raisebox{0.0\height}{\includegraphics[width=0.6\textwidth]{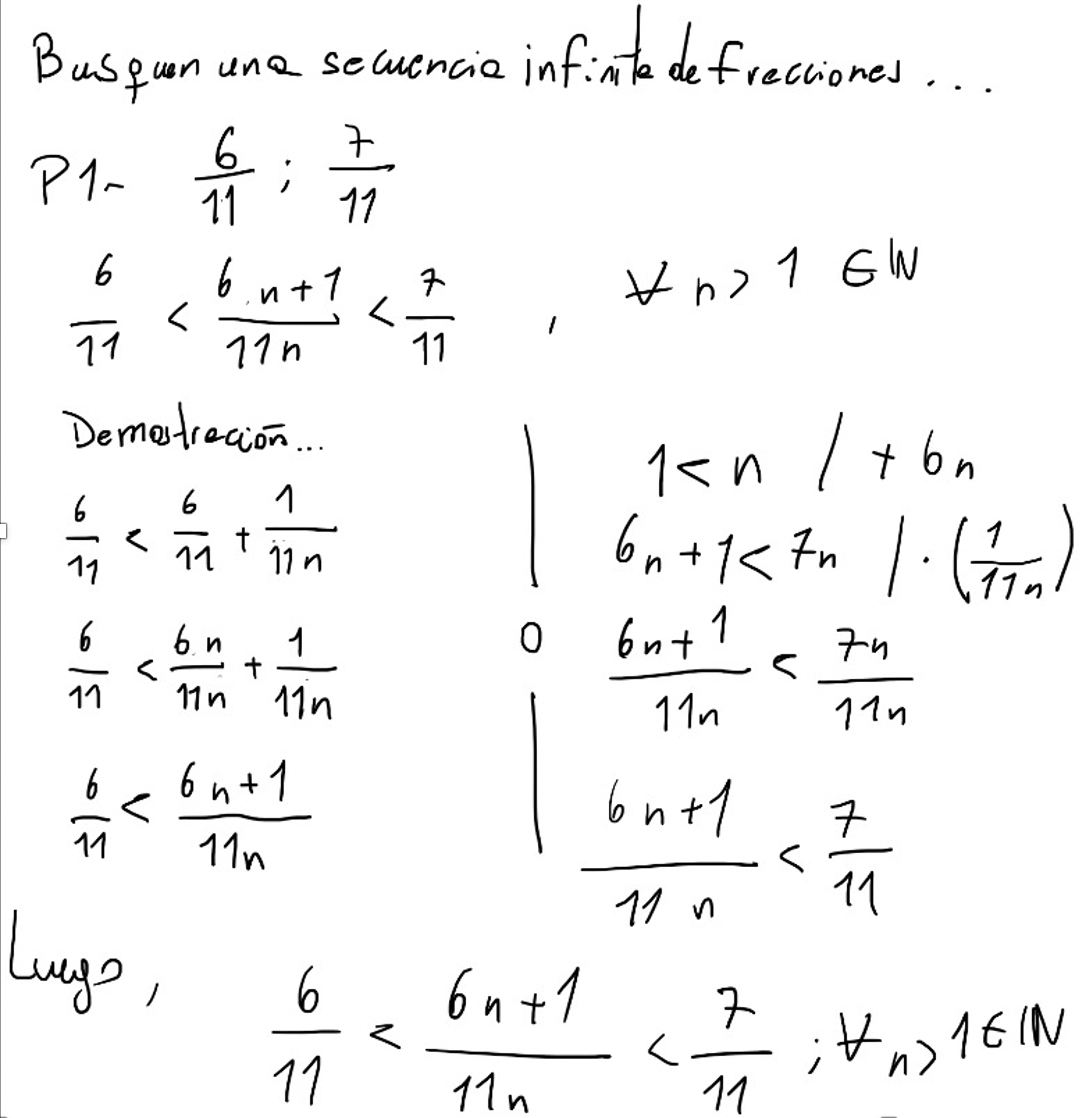}}
\caption{Una de las demostraciones enviadas por una estudiante para
	justificar que la secuencia de infinitos valores est{\'a} entre
	$\dfrac{6}{11}$ y $\dfrac{7}{11}$.}
	\end{center}

\end{figure}

\subsubsection{Justificaciones del segundo grupo de secuencias}

El {\'u}nico grupo de respuestas que utiliz{\'o} argumentos de naturaleza
diferente fue el grupo 2 (entre D6 y D8 de la tabla de la Figura 11). En
efecto, utilizaron Geogebra para graficar la funci{\'o}n continua asociada y a
partir de la gr{\'a}fica justificar que se encuentra acotada. Aunque la
justificaci{\'o}n es escueta, se puede apreciar en ella un trabajo en el plano
semi{\'o}tico-instrumental. El artefacto digital --y particularmente su
sistema de gr{\'a}ficas-- funciona como elemento principal de la
justificaci{\'o}n.

\begin{figure}[h!]
	\begin{center}
  \raisebox{0.0\height}{\includegraphics[width=0.7\textwidth]{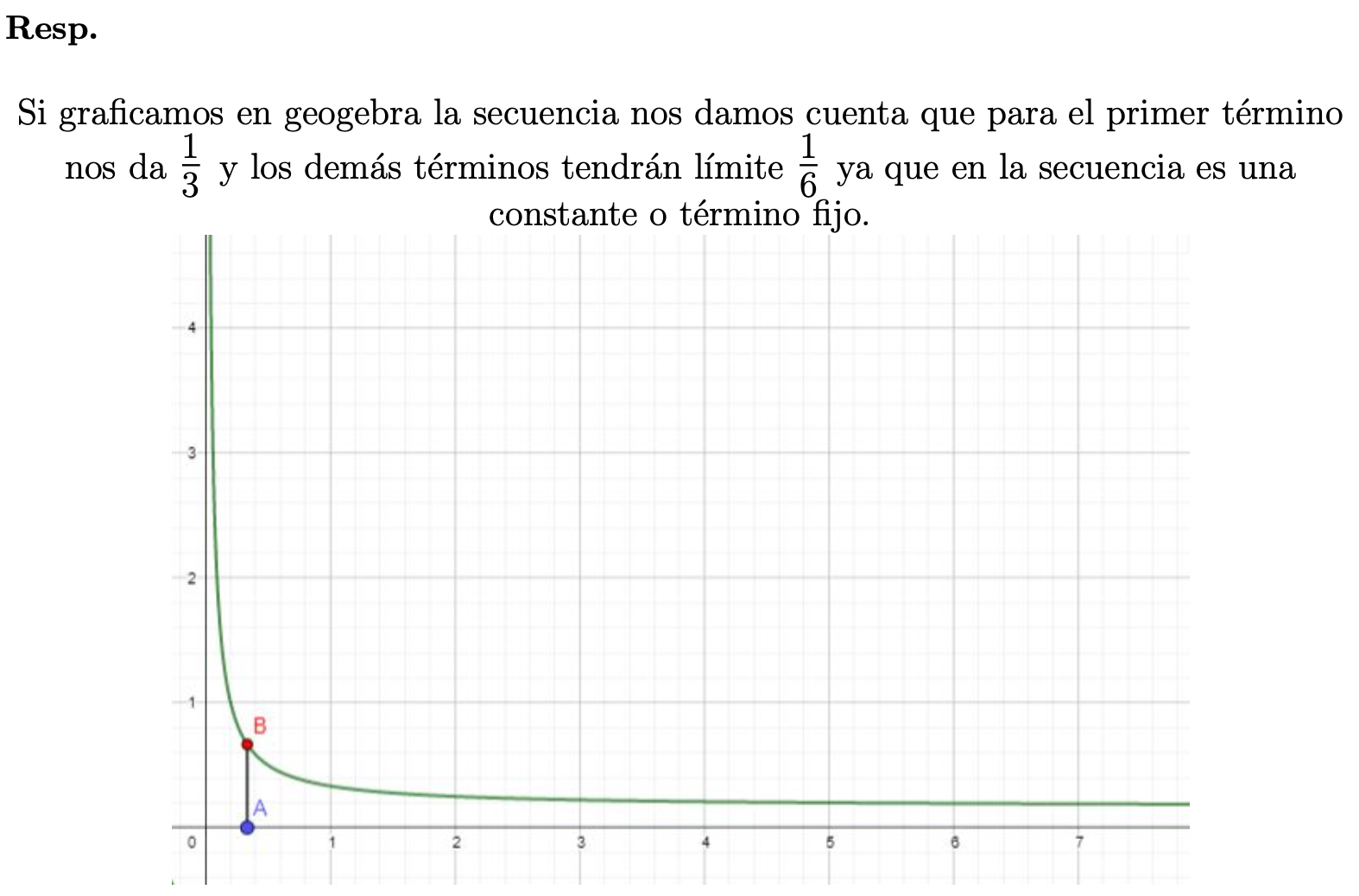}}
  \caption{Justificaci{\'o}n de la respuesta D6, donde se utiliza Geogebra
  para graficar la funci{\'o}n $\dfrac{6}{8} - \dfrac{1}{8 x}$}
\end{center}
\end{figure}

El esquema de la Figura 12 muestra la transformaci{\'o}n del proceso
discursivo. Se puede observar que en un comienzo el uso coordinado de
artefactos simb{\'o}licos y elementos del referencial te{\'o}rico para
construir una justificaci{\'o}n, provoca la activaci{\'o}n principalmente del
plano vertical instrumental-discursivo, salvo una respuesta que activa el
plano semi{\'o}tico-instrumental. Posteriormente se construyen demostraciones,
lo que evidencia un enriquecimiento de la g{\'e}nesis discursiva en
relaci{\'o}n con lo que se observ{\'o} en la primera clase. Finalmente, los
argumentos de las demostraciones son instrumentalizados, pues son utilizados
posteriormente, por casi todos los estudiantes, para dar respuesta a nuevas
tareas de una forma similar.

Si resumimos la activaci{\'o}n de las distintas g{\'e}nesis y planos, durante
la clase 2 y la evaluaci{\'o}n, se puede apreciar un tr{\'a}nsito desde lo
discursivo hacia lo instrumental; las justificaciones, en su mayor{\'i}a usan
un esquema argumentativo similar, salvo el grupo 3 que utiliza una
argumentaci{\'o}n usando el sistema de gr{\'a}ficas de Geogebra.

\begin{figure}[h!]
  \raisebox{0.0\height}{\includegraphics[width=0.99\textwidth]{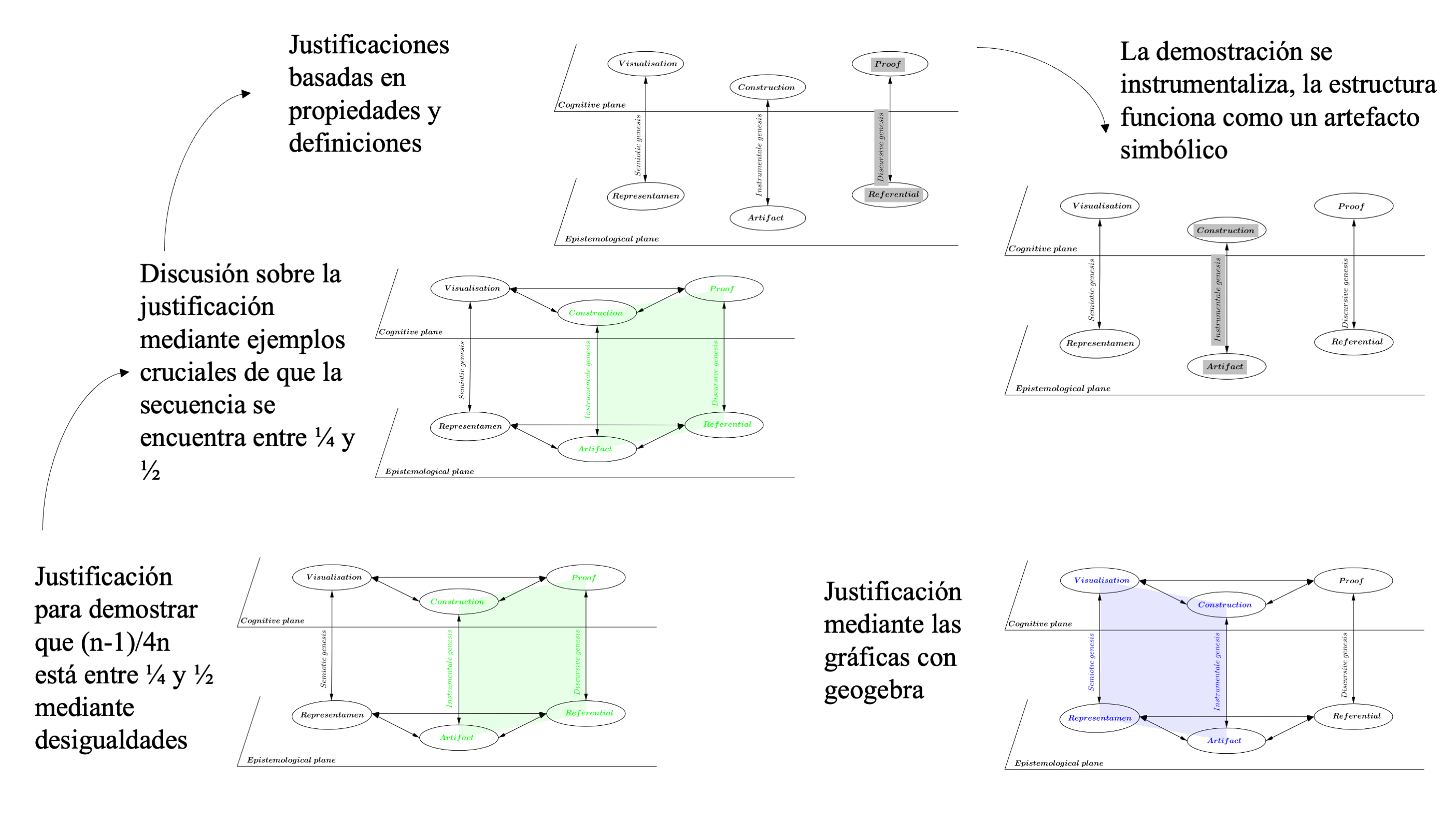}}
  \caption{Esquema del trabajo matem{\'a}tico en la clase 2 y de la
  evaluaci{\'o}n sobre la construcci{\'o}n de secuencias infinitas que
  est{\'e}n entre 1/4 y 2/4 y luego entre dos fracciones definidas
  aleatoriamente por la plataforma.}
\end{figure}

\section{Conclusi{\'o}n}
En las clases virtuales realizadas durante el a{\~n}o 2020 se observ{\'o} poca
participaci{\'o}n de los estudiantes en las clases; las intervenciones eran
escasas, y la mayor{\'i}a de las veces las c{\'a}maras y micr{\'o}fonos se
mantienen apagados. Entre las iniciativas adoptadas por docentes para la
formaci{\'o}n de futuros profesores, se comenzaron a desarrollar preguntas
abiertas, que tuvieran infinitas soluciones, con correcci{\'o}n y feedback
autom{\'a}tico, con el objetivo de fomentar la discusi{\'o}n sobre las
estrategias m{\'a}s que sobre las respuestas. Los resultados del trabajo
matem{\'a}tico de los estudiantes, luego de implementar una secuencia de
tareas sobre fracciones, son los que se presentan en este art{\'i}culo.

Para el estudio nos planteamos abordar la pregunta:
{\textquestiondown}cu{\'a}l es el trabajo de argumentaci{\'o}n puesto en juego
en un contexto de clase virtual donde se usa una tarea en un artefacto de
evaluaci{\'o}n en l{\'i}nea para una tarea de fracciones?, y para poder
responderla, se realiz{\'o} una intervenci{\'o}n compuesta por dos clases y
una evaluaci{\'o}n.

\

\

\

\

En la intervenci{\'o}n, se les pidi{\'o} primero a los estudiantes encontrar
una fracci{\'o}n de la forma $\dfrac{a}{b}$ con $a$ y $b$ en los enteros,
entre $\dfrac{1}{4}$ y $\dfrac{2}{4}$. Por la densidad de los racionales, esta
pregunta tiene infinitas soluciones y una vez que se ingresaba la respuesta,
el sistema evaluaba si se encontraba o no en el intervalo y si cumpl{\'i}a o
no con el formato pedido, y entregaba un feedback. Adem{\'a}s, daba razones de
por qu{\'e} la respuesta era correcta. La situaci{\'o}n evolucion{\'o} a
partir del trabajo de discusi{\'o}n y finaliz{\'o} con una pregunta que
ped{\'i}a encontrar una secuencia con infinitas fracciones, primero entre \
$\dfrac{1}{4}$ y $\dfrac{2}{4}$ y luego entre dos fracciones con valores
aleatorios dados por la plataforma.

Luego del trabajo individual en la plataforma, se comenz{\'o} a discutir con
los estudiantes sobre las estrategias utilizadas, y estos solo declararon dos
estrategias: con decimales, y amplificando por el mismo n{\'u}mero el
numerador y denominador de cada fracci{\'o}n. Esta {\'u}ltima estrategia hizo
aparecer la siguiente conjetura: si $\dfrac{1}{4}$ y $\dfrac{2}{4}$ se
amplificaba por n entonces aparec{\'i}an $n - 1$ fracciones entre ellas. Uno
de los estudiantes hizo la demostraci{\'o}n al finalizar la clase.

Adem{\'a}s, en la discusi{\'o}n apareci{\'o} otra estrategia que fue por
ensayo y error, que permiti{\'o} entrar en una discusi{\'o}n epist{\'e}mica.
Esta discusi{\'o}n se caracteriz{\'o} por un tr{\'a}nsito desde argumentos
instrumentales --tales como ``elevando la fracci{\'o}n a menos 1 se da vuelta
la desigualdad'', argumentos err{\'o}neos, y algunos argumentos circulares--
hacia la g{\'e}nesis discursiva (Kuzniak et al., 2016a), basados en
propiedades del referencial te{\'o}rico.

En una segunda clase se extendi{\'o} la pregunta hacia la b{\'u}squeda de
secuencias infinitas de fracciones que estuvieran entre $\dfrac{1}{4}$ y
$\dfrac{2}{4}$. Esta tarea fue m{\'a}s compleja para los estudiantes, pero
hubo al menos dos grupos que encontraron algunas secuencias correctas. Se
discuti{\'o} sobre la diferencia entre mostrar que una secuencia cumple para
varios valores y demostrar para todo n natural, lo que Stylianides y
Stylianides (2009, p. 315) llaman argumentos emp{\'i}ricos versus una
demostraci{\'o}n y que, dependiendo de la tarea y de c{\'o}mo se implemente,
permite a los estudiantes construir un proceso de justificaci{\'o}n en el cual
busquen argumentos en un referencial te{\'o}rico.

Finalmente, en la evaluaci{\'o}n se pidi{\'o} a cada estudiante ingresar en la
plataforma un valor entre fracciones aleatorias, digamos $\dfrac{a}{b}$ y
$\dfrac{c}{d}$. Luego, se pide encontrar una secuencia infinita que estuviera
entre $\dfrac{a}{b}$ y $\dfrac{c}{d}$, y demostrar que tal secuencia estaba
entre ambas fracciones dadas. La respuesta y la correcci{\'o}n de las
{\'u}ltimas dos tareas se realizaba en forma manual. El argumento principal
para demostrar que una secuencia infinita estaba entre dos fracciones fue
an{\'a}logo al mostrado en la clase 1, dicho de otra forma, la estructura de
la argumentaci{\'o}n pas{\'o} a ser un artefacto simb{\'o}lico para todos los
grupos.

Se advierte que es dif{\'i}cil saber qu{\'e} camino tomar{\'a}n estas
discusiones y que dependen en gran medida del papel mediador del profesor. Tal
como lo indica Yackel (2002), es rol del profesor reconocer la importancia y
validez de los argumentos, conocer las posibilidades conceptuales de los
estudiantes y conocer los conceptos matem{\'a}ticos subyacentes para transitar
hacia una argumentaci{\'o}n matem{\'a}tica significativa. Es dif{\'i}cil
pensar en una planificaci{\'o}n a priori para saber hacia d{\'o}nde
llevar{\'a} la discusi{\'o}n por lo que la posibilidad de replicar la
intervenci{\'o}n no es lo que se demuestra, sino el proceso de
argumentaci{\'o}n en s{\'i} mismo.

Durante toda la intervenci{\'o}n se percibe un juego interesante entre las
g{\'e}nesis instrumentales y discursivas. Si en la clase 2 aparecieron
argumentos en la dimensi{\'o}n discursiva, estos mismos argumentos pasaron a
ser artefactos simb{\'o}licos cuando justificaron sus respuestas en la
evaluaci{\'o}n. Como la mayor{\'i}a de los estudiantes utilizaron el mismo
esquema argumentativo, se podr{\'i}a interpretar como ``una secuencia
culturalmente codificada de acciones que se instancian continuamente en la
pr{\'a}ctica social'' (Radford, 2014, p. 417). Los argumentos construidos por
un estudiante en la clase 2 fueron compartidos con los compa{\~n}eros, quienes
adaptaron la justificaci{\'o}n a cada par de n{\'u}meros que les asign{\'o} la
plataforma de forma aleatoria. Si bien se pudo observar que en diferentes
instancias se activa el plano instrumental-discursivo (Kuzniak y Richard,
2014), las cualidades de los argumentos sufrieron transformaciones gracias a
la utilizaci{\'o}n de diferentes artefactos y elementos del referencial
te{\'o}rico.

\

\

El rol de la tecnolog{\'i}a es devolver a los estudiantes la responsabilidad
de comprometerse con una respuesta, de tal forma que cada uno tenga algo que
decir. En clases previas, de alguna forma, los estudiantes no se sent{\'i}an
comprometidos a responder a las preguntas planteadas, las discusiones eran
escasas y forzadas. Seg{\'u}n Solar y Piquet (2016) una de las condiciones,
para promover la argumentaci{\'o}n, es dar oportunidades de participaci{\'o}n.
En este sentido, la plataforma permite dar estas oportunidades de forma
expl{\'i}cita. Tambi{\'e}n, cuando la plataforma interviene con correcci{\'o}n
y feedback autom{\'a}tico, la informaci{\'o}n que entrega al sujeto genera
significados, los cuales hay que cuestionar, de manera de seleccionar aquellos
que son epistemol{\'o}gicamente v{\'a}lidos o, si no lo son, tener clara su
validez relativa.

La literatura analizada muestra que es necesario trabajar en el desarrollo de
la argumentaci{\'o}n de futuros profesores, tanto porque es una habilidad de
orden superior, como porque los trabajos en profesores debutantes muestran
debilidades al respecto. Adem{\'a}s, gracias a que se est{\'a}n formando como
profesores de matem{\'a}ticas, se advierte c{\'o}mo el contrato did{\'a}ctico
(Brousseau, 1998) los obliga a buscar argumentos basados en propiedades
matem{\'a}ticas m{\'a}s que a solo resolver, por lo que es posible que esta
misma situaci{\'o}n funcione en estudiantes de otros niveles educativos.

Una de las perspectivas de este trabajo es ver qu{\'e} otras discusiones
emergen en distintos contextos formativos y si es posible automatizar otras
preguntas que, eventualmente, puedan aparecer. Por ejemplo, que un sistema sea
capaz de validar si una secuencia est{\'a} acotada, en particular, que se
pueda automatizar la validaci{\'o}n de que una secuencia \{an\} est{\'a}
acotada por dos fracciones, esto para que sirva como forma de
experimentaci{\'o}n para los estudiantes, para reducir los tiempos de
correcci{\'o}n de los profesores, y a trav{\'e}s del feedback, instar a la
b{\'u}squeda de secuencias con otras caracter{\'i}sticas distintas a las
observadas en esta investigaci{\'o}n.

\end{document}